%% file: main.tex
\documentclass[12pt, a4paper]{article}

\usepackage[T1]{fontenc}
\usepackage[utf8]{inputenc}
\usepackage{lmodern}
\usepackage[english]{babel}

\usepackage[margin=3cm]{geometry}
\usepackage[x11names]{xcolor}

\usepackage{adjustbox}
\usepackage{amsmath,amsthm,amssymb,mathrsfs}
\usepackage{array, makecell}
\usepackage{authblk}
\usepackage{bbm}
\usepackage{blindtext}
\usepackage{bm}
\usepackage{cancel}
\usepackage{censor}
\usepackage{collcell}
\usepackage{comment}
\usepackage{enumitem}
\usepackage{epigraph}
\usepackage{esvect}
\usepackage{etoolbox}
\usepackage{fancyhdr}
\usepackage[bottom]{footmisc}
\usepackage{graphicx}
\usepackage[pdfborder={0 0 0}, bookmarksdepth=-1, unicode=true]{hyperref}
\usepackage{indentfirst}
\usepackage{index} 
\usepackage{lipsum}
\usepackage{mathtools}
\usepackage{multicol}
\usepackage[square, comma, numbers, super, sort&compress]{natbib}
\usepackage{nicematrix}
\usepackage{pdfpages}
\usepackage{pgfplots}
\usepackage{pifont}
\usepackage{ragged2e}
\usepackage{scalerel}
\usepackage{setspace}
\usepackage{tabularx}
\usepackage{tcolorbox}
\usepackage{thmtools}
\usepackage{tikz}
\usepackage{titletoc}
\usepackage{type1cm}
\usepackage[normalem]{ulem}
\usepackage{wrapfig}
\usepackage{xparse}
\usepackage{xstring}
\usepackage{xurl}
\usepackage{yhmath}

\pgfplotsset{compat=1.15}
\usetikzlibrary{shapes.geometric}
\usetikzlibrary{calc,shadows,patterns,angles,quotes}
\usetikzlibrary{arrows}

\newcommand{\resetcounters}
  {
  \setcounter{section}{0}
  \setcounter{ans}{0}
  \setcounter{ax}{0}
  \setcounter{con}{0}
  \setcounter{clm}{0}
  \setcounter{ex}{0}
  \setcounter{pr}{0}
  \setcounter{pp}{0}
  \setcounter{prop}{0}
}

\newtheoremstyle{one}
  {6pt}{15pt}{}{}{\bfseries}{.}{1em}{}
\theoremstyle{one}

\newtheorem{con}{Conclusion}[section]

\newtheorem{pr}{Problem}[subsection]

\newtheorem{rem}{Remark}[section]
\newtheorem{ex}{Example}[section]

\newtheorem{pp}{Property}[subsection]

\newtheorem{df}[pr]{Definition}

\newtheorem{cj}{Conjecture}

\newtheorem{cons}[pr]{Construction}
\newtheorem{alg}[pr]{Algorithm}

\newenvironment{pf}[1][0]
 {%
  \begin{proof}[\textit{Proof\ifnum#1>0 \, #1\fi.}]
 }{\end{proof}}

\newsavebox\curwrapfig
\makeatletter
\long\def\wrapfiguresafe#1#2#3{%
  \sbox\curwrapfig{#3}%
  \par\penalty-100%
  \begingroup
    \dimen@\pagegoal \advance\dimen@-\pagetotal
    \advance\dimen@-\baselineskip
    \ifdim \ht\curwrapfig>\dimen@
      \break%
    \fi%
  \endgroup%
  \begin{wrapfigure}{#1}{#2}%
    \usebox\curwrapfig%
  \end{wrapfigure}%
}
\makeatother

\title{Constructions of Large m-Distance Sets on Triangular Lattice}
\author[1]{Li-Ren Bao}
\author[2]{Wei-Hsuan Yu}
\affil[1]{Taipei Municipal Chien Kuo High School, Taipei, Taiwan. \protect\\Email: richman960118@gmail.com}
\affil[2]{Department of Mathematics, National Central University, Taoyuan, Taiwan.\protect\\
\small\itshape{Email: whyu@math.ncu.edu.tw}}

\begin{document}
\lhead{}
\chead{}
\rhead{}
\cfoot{\thepage}
\rfoot{}

\input{name.tex}

\bibliographystyle{plain}
\bibliography{ref}

\end{document}

%% file: name.tex
\date{\today}
\lfoot{}
\resetcounters
\maketitle

\setlist[enumerate]{topsep=0pt}
\setlist[itemize]{topsep=0pt}

\pagenumbering{Roman}
\setcounter{page}{1}
\pagenumbering{arabic}
\setcounter{page}{1}

\begin{abstract}
    An $m$-distance set is a collection of points such that the distances between any two points have $m$ possible values. We use two different methods to construct large $m$-distance sets on the triangular lattices. One is to use the first m smallest distances and find the largest cliques, and the other is using the notions of hexagons. Multiplicities of the distances were observed for comparison for the two methods.
\end{abstract}

\section{Introduction}

An $m$-distance set is a point set where any two points have only $m$ possible distances between them. If an $m$-distance set gains an $m+1$-th distance when a new point is added, we will call it a maximal $m$-distance set. If an $m$-distance set has the largest number of points among all possible $m$-distance sets, it is called a maximum $m$-distance set. Note that a maximal $m$-distance is not always maximum, but a maximum $m$-distance set is always maximal.

Determining the size of the maximum $m$-distance set in $\mathbb R^n$ has been difficult, and we only know values for small $m$. Erdős and Fishburn \cite{erdHos1996maximum} determined that the maximum sizes of planar $2, 3, 4$, and $5$ distance sets are $5, 7, 9$, and $12$ points, respectively. They also found all possible constructions for $2,3,4$ distance sets in the paper. Then, Shinohara finished the classification of planar $3$-distance sets \cite{shinohara2004classification}, and then proved the maximum planar $5$-distance set is unique \cite{shinohara2008uniqueness}. Finally, Wei proved that the maximum planar $6$-distance set has size $13$ \cite{xianglin2012proof}.

For higher dimensions, to determine the maximum $m$-distance set is much more difficult. However, for $2$-distance sets, Lisoněk determined the maximum size in $\mathbb R^d$ for $d\leq 8$ \cite{lisonvek1997new}. Bannai, Sato, Shigezumi also gave a construction method for maximal $m$-distance sets in $\mathbb R^n$ based on Johnson Graphs \cite{bannai2012maximal}. Other than Johnson Graphs, Adachi, Hayashi, Nozaki, and Yamamoto also classified some maximal $m$-distances set containing the Hamming Graph \cite{adachi2017maximal}.

Here, we will discuss the construction of points located on the triangular lattices in a plane. These are points with coordinates $a(1, 0)+b(\frac12, \frac{\sqrt3}2)$, where $a, b\in\mathbb Z$. We discuss the problem in triangular coordinates due to the following conjecture by Erdős and Fishburn \cite{erdHos1996maximum}.
\begin{cj}
    For any $m\geq 3$, at least one of the maximum $m$-distance sets lies within triangular lattices. If $m\geq7$, then maximum $m$-distance sets can only be found in triangular lattices.
\end{cj}
For the cases of $m \leq 6$, the above conjecture is true. 

There are some recent studies on the cases of triangular lattices. Most of them are based on special hexagonal figures. Balaji, Edwards, Loftin, Mcharo, Phillips, Rice, and Tsegaye gave a bound with regular polygons and especially regular hexagons in triangular lattices \cite{balaji2019lattice}. Ahmed and Wildstrom used some hexagons that are symmetrical with respect to vertical lines and bounded the number of distances inside them \cite{ahmeddistance}.


Here, we determine the maximum $m$- distance sets from $m=7$ to $34$ under the assumption of using the $m$ smallest distances in the triangular lattices. 
We also present an alternate method of construction by equiangular hexagons. If we use equiangular hexagons instead of regular hexagons, we can find constructions of larger size than the case found in \cite{balaji2019lattice}. The following values are the best lower bounds we found for sets that determine at most $m$ distances by  above two methods. The stars mean the construction method of hexagons is giving strictly larger constructions.

\begin{center}
\begin{tabular}{|c|c|c|c|c|c|c|c|c|c|c|c|c|c|c|}
\hline
$m$ & $7$ & $8$ & $9$ & $10$ & $11$ & $12$ & $13$ & $14$ & $15$ & $16$ & $17$ & $18$ & $19$ & $20$\\\hline
size & $16$ & $19$ & $21$ & $24$ & $27$ & $27$ & $31$ & $34$ & $37$ & $37$ & $42$ & $45$ & $48$* & $49$\\\hline
$m$ & $21$ & $22$ & $23$ & $24$ & $25$ & $26$ & $27$ & $28$ & $29$ & $30$ & $31$ & $32$ & $33$ & $34$\\\hline
size & $55$ & $58$ & $61$* & $63$ & $63$ & $69$ & $72$ & $75$* & $75$* & $79$ & $79$ & $85$ & $88$ & $91$*\\\hline
\end{tabular}
\end{center}

When designing the first method, we assume that the smallest $m$ distances in triangular lattices will yield the maximum $m$-distance set. Therefore, we determine all points that have those distances with both $(0, 0)$ and $(0, 1)$. Then we build a graph by using these distances to determine the edges, and find the maximum clique to construct the $m$-distance set. 

As these are in fact maximum when we are discussing distance sets based on the $m$ smallest distances, we have the following conjecture
\begin{cj}
    The maximum $7$-distance and $8$-distance sets have $16$ and $19$ points, respectively. The $7$-distance set can only be realized by taking all triangular lattices in a regular hexagon with side length $2$, and removing three all adjacent or all non-adjacent corners. The $8$-distance set can only be realized by taking all points in that regular hexagon.
\end{cj}
However, when we construct the $19$-distance set, we noticed that the smallest 19 distances cannot form a maximum size $19$-distance set. The $19$ smallest distances would only yield a $19$-distance set of $45$ points, while an equiangular hexagon with side length alternating $3$ and $4$ will form a $19$-distance set of $48$ points. The main reason is that when we construct the $19$-distance set by the $19$ smallest distances, the $19$th smallest distance only has multiplicity $2$. The hexagon eliminated this distance and added the $20$th smallest one which is $7$, resulting in a construction with $3$ more points. This observation gives two other conjectures. 
\begin{cj}
    Taking all points in an equiangular hexagon with side length alternating $k$ and $k+1$, or a regular hexagon, will always give a maximum m-distance set.
\end{cj}
There are two confirmed examples for this conjecture: All triangular lattice points in a regular hexagon with side length $1$ give a maximum $3$-distance set, and all triangular lattice points in an equiangular hexagon alternating side lengths $1$ and $2$ give the unique maximum $5$-distance set.
\begin{cj}
    For any $m\geq 3$, at least one construction of maximum $m$-distance set can be realized by repeatedly removing points from the aforementioned hexagons.
\end{cj}
This can be confirmed with a maximum $4$-distance set construction, which can be found by removing $3$ corners from an equiangular hexagon with side length alternating $1$ and $2$, and also with a maximum $6$-distance set, which can be obtained by deleting $6$ corners from a regular hexagon with side length $2$.
In the end, we found that for $m=19,23,28,29,34$, the construction with equiangular hexagons has more points than the construction with the maximum clique method. The main reason is that some distances have low multiplicities if we use the $m$ smallest distances in the maximum clique method. Substituting distances of small multiplicities out for a bigger integer yields a construction of larger size by the hexagon. This gives yet another conjecture about the multiplicities.

\begin{cj}
    A maximum $m$-distance set will not have a distance with multiplicity less than or equal to $2$.
\end{cj}

\section{Preliminary}
In this paper, we use the following method to check the size of a maximum $m$-distance set, assuming the $m$ distances are given. The method is based on the paper of Szöllősi and Östergård
\cite{szollHosi2018constructions}.

\begin{alg}
    Given a set $\mathcal S$ with $m$ distances, we will consider the two points $A(0, 0), B(0, 1)$, and find all points $P$ such that $d(P, A)\in \mathcal S$ and $d(P, B)\in \mathcal S$. We would put these points in a set called $V$. Now we build a graph with the vertex set $V$, and two vertices share an edge if the distance between the two points is in $\mathcal S$. Finally, by finding the largest clique $W$, we have $W\cup\{A, B\}$ being an $m$-distance set.
\end{alg}

\begin{pf}
    We know that any two points in the largest clique have a distance in $\mathcal S$, and any point in the largest clique also has a distance to $A$ and $B$ in $\mathcal S$. Since $|\mathcal S|=m$ we know that any two points among the construction have at most $m$ possible distances.
\end{pf}
By setting the distances in $\mathcal S$ to the distances in the triangular lattice, we can find the largest $m$-distance set by building the graph and running a maximum clique algorithm on the graph. The maximum graph algorithm we used here is designed by Carraghan and Pardalos
\cite{carraghan1990exact}.

\begin{ex}
    If we want to find the largest $7$-distance set by smallest $7$ distances in of triangular lattice. We will first take $A(0, 0), B(0, 1)$, and then locate all points $P$ with $d(P, A)\in\{1, \sqrt3, 2, \sqrt7, 3, \sqrt{12}, \sqrt{13}\}$ and $d(P, B)\in\{1, \sqrt3, 2, \sqrt7, 3, \sqrt{12}, \sqrt{13}\}$ and take them into a set $V$. For example, $(3, 0)$ would be in $V$ as it has distance $3$ with $A$ and $2$ with $B$. Then we will add an edge between two points if they have a distance of one of $1,\sqrt3, 2, \sqrt7, 3, \sqrt{12}, \sqrt{13}$ with each other. For example, $(2, 0)$ and $(\frac12, \frac{\sqrt3}2)$ would have an edge between them because they have a distance of $\sqrt3$.

    Finally, we will find the maximum clique in such graph. These points would always have distances one of $1,\sqrt3, 2, \sqrt7, 3, \sqrt{12}, \sqrt{13}$ with each other, with $A$, and with $B$. Therefore, the maximum clique plus $A$ and $B$ would form the largest $7$-distance set with the aforementioned distances.
\end{ex}
We implemented the algorithms and got the results by the codes in the following sections. All codes are available at https://github.com/BaoCoder613/triangular-m-distance. The sequence of the distances of set $\mathcal S$ in the code is taken from OEIS \cite{oeis:A003136}.

\section{Constructions With the Smallest Distances}
\subsection{$m=7$ to $16$}
Here we present the construction of $7$-distance to $16$-distance set, with the $m$ smallest triangular lattice distances.
Some images are omitted. Those are the cases where the construction of $m$-distance set does not contain more points than $m-1$-distance set, and the code outputs are the exact same construction for them.
Some discussions, especially geometrical properties and relations with other distance sets on the constructions, are listed below:

\begin{cons}[$7$-distance, $16$ points]
    The construction is three corner points off of a hexagon with side length $2$, giving it $3$ axes of symmetry. It can also be formed by adding $3$ points from the known $6$-distance construction, or deleting $3$ corners from the $8$-distance construction.
\end{cons}
\begin{cons}[$8$-distance, $19$ points]
    The construction is a regular hexagon with side length $2$, giving it $6$ axes of symmetry.
\end{cons}
\begin{cons}[$9$-distance, $21$ points]
    Using this method of finding the largest clique, we found that it does not contain the previous construction of $m=8$. However, it can be formed by deleting $6$ corners from the $11$-distance. This can be further discussed by considering all the possible maximum $9$-distance sets in triangular lattices. It still has line symmetry.
\end{cons}
\begin{center}
\begin{tabular}{ccc}
    \includegraphics[width=0.3\linewidth]{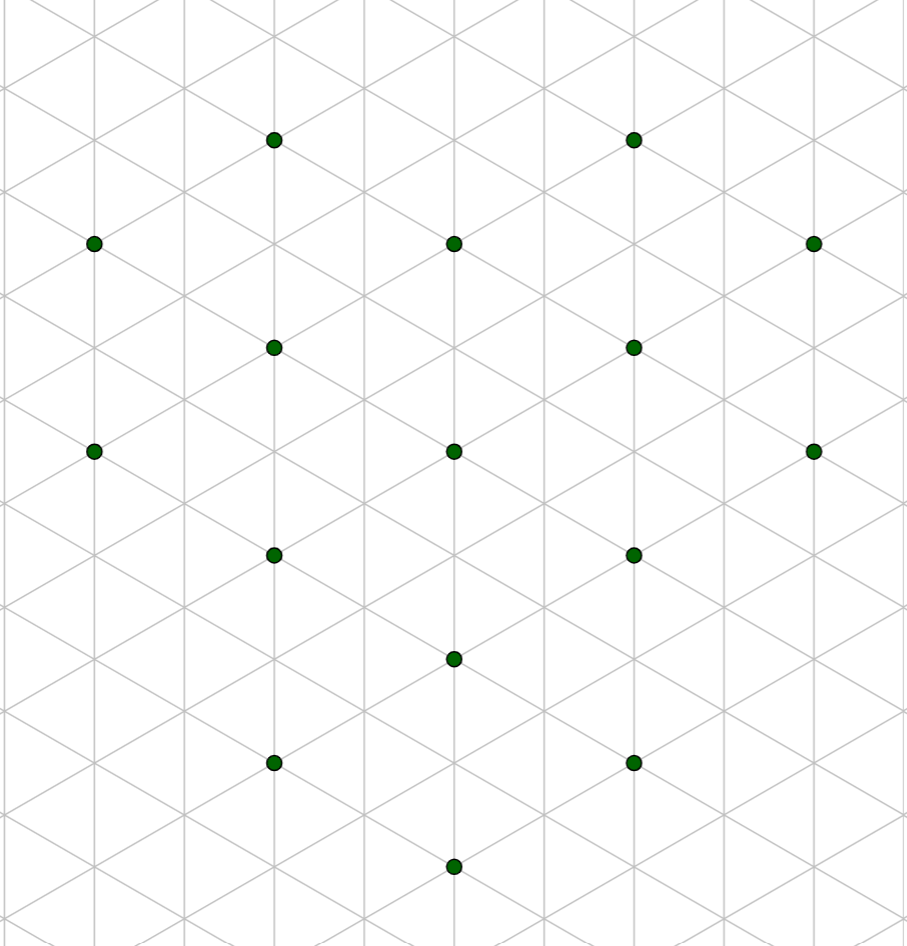}&
\includegraphics[width=0.3\linewidth]{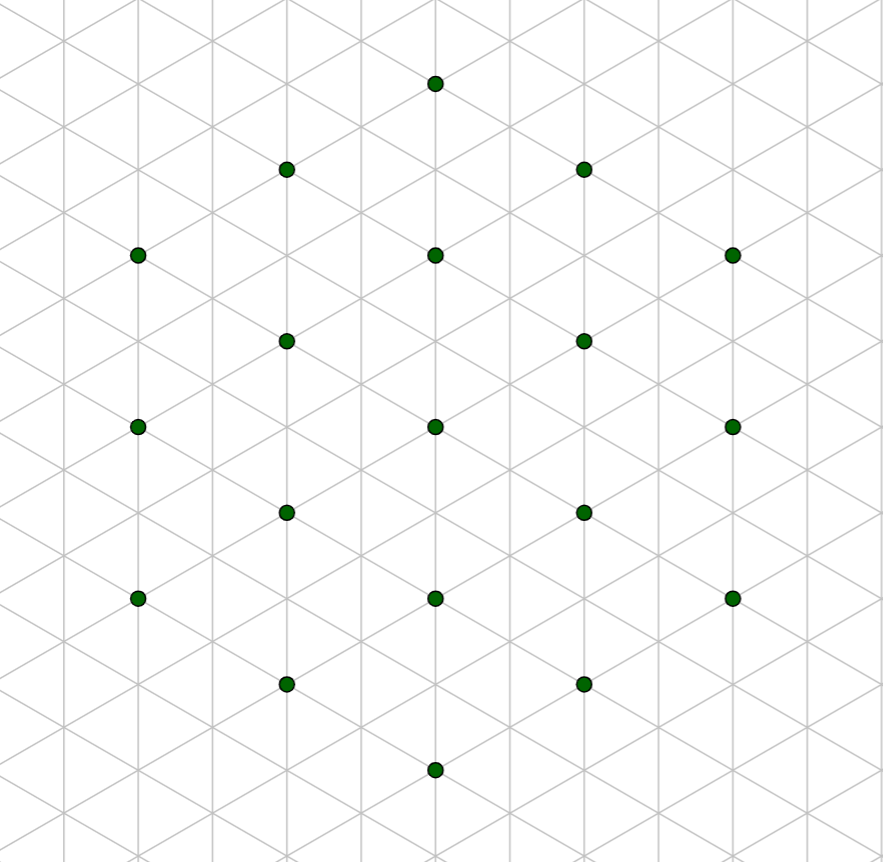}&
\includegraphics[width=0.3\linewidth]{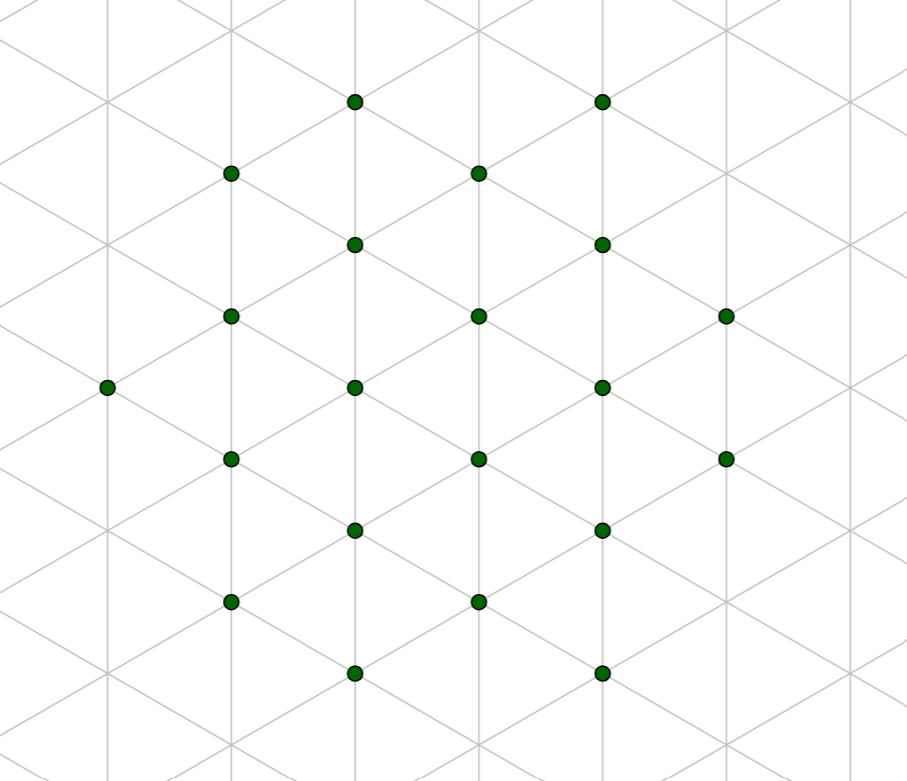}\\
$7$-distance & $8$-distance & $9$-distance
\end{tabular}
\end{center}

\begin{cons}[$10$-distance, $24$ points]
    The construction is three corner points off of an equiangular hexagon with side length $2$ and $3$ alternating. This one is asymmetrical.
\end{cons}
\begin{cons}[$11$-distance, $27$ points]
    This figure is an equiangular hexagon, with three sides of length $2$ while the others are of length $3$, and the lengths are alternating. It also has $3$ axes of symmetry.
\end{cons}
\begin{cons}[$12$-distance, $27$ points]
    The code got the exact same construction with $11$-distance.
\end{cons}
\begin{cons}[$13$-distance, $31$ points]
    It is a configuration that removes $6$ corners from a regular hexagon with side length $3$, or the $15$-distance construction. which suggests a similar structure to $6$-distance. It has $6$ axes of symmetry.
\end{cons}

\begin{center}
\begin{tabular}{ccc}
\includegraphics[width=0.3\linewidth]{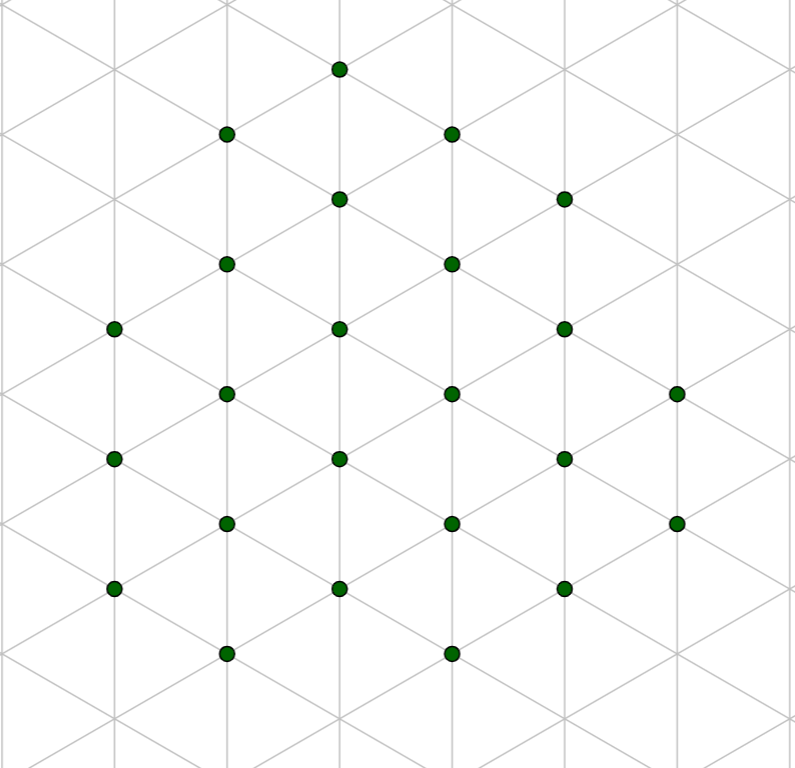}&
\includegraphics[width=0.3\linewidth]{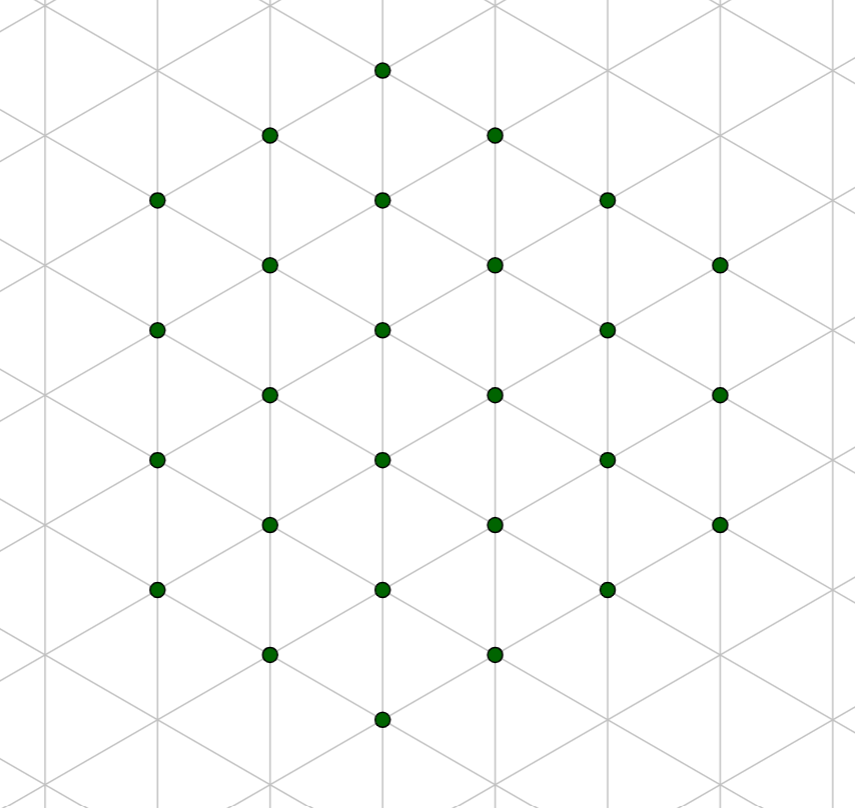}&
\includegraphics[width=0.3\linewidth]{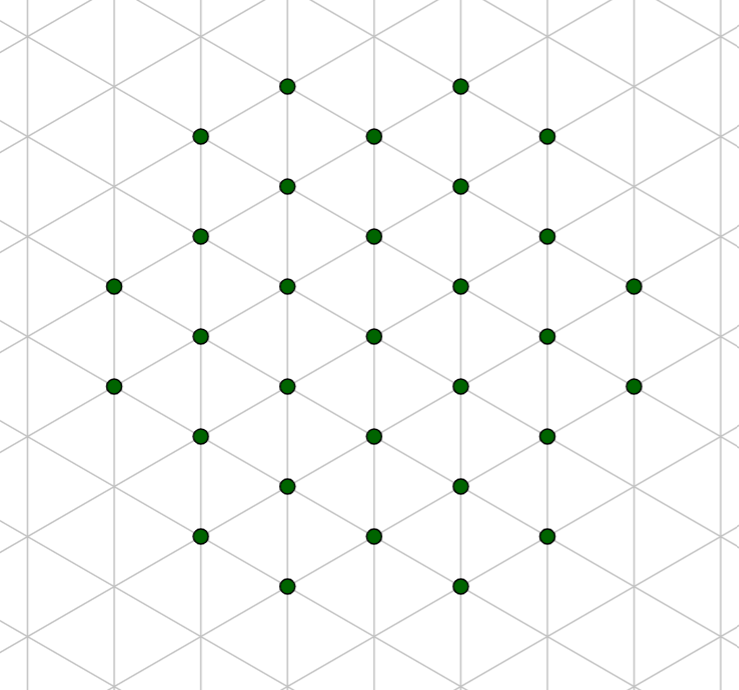}\\
$10$-distance & $11$-distance & $13$-distance
\end{tabular}
\end{center}

\begin{cons}[$14$-distance, $34$ points]
    The figure is removing $3$ corners from a regular hexagon with side length $3$. Unlike $7$-distance, though, the removed corners are adjacent.
\end{cons}
\begin{cons}[$15$-distance, $37$ points]
    This is a regular hexagon with side length $3$ and $6$ axes of symmetry.
\end{cons}
\begin{cons}[$16$-distance, $37$ points]
    The construction is the same as $m=15$ again. This can suggest that a highly symmetrical hexagonal construction with the $m$-smallest distance may lead to the maximum $m+1$-distance set to be the same size as the maximum $m$-distance set.  When we add one more point to this hexagon, it will create two new distances. 
\end{cons}

\begin{center}
\begin{tabular}{cc}
\includegraphics[width=0.3\linewidth]{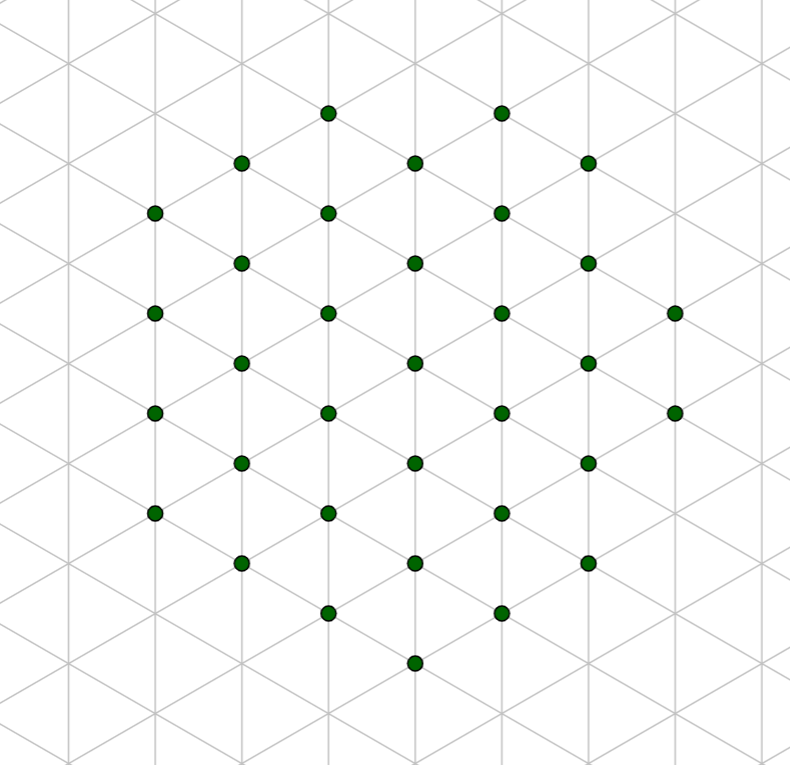}&
\includegraphics[width=0.3\linewidth]{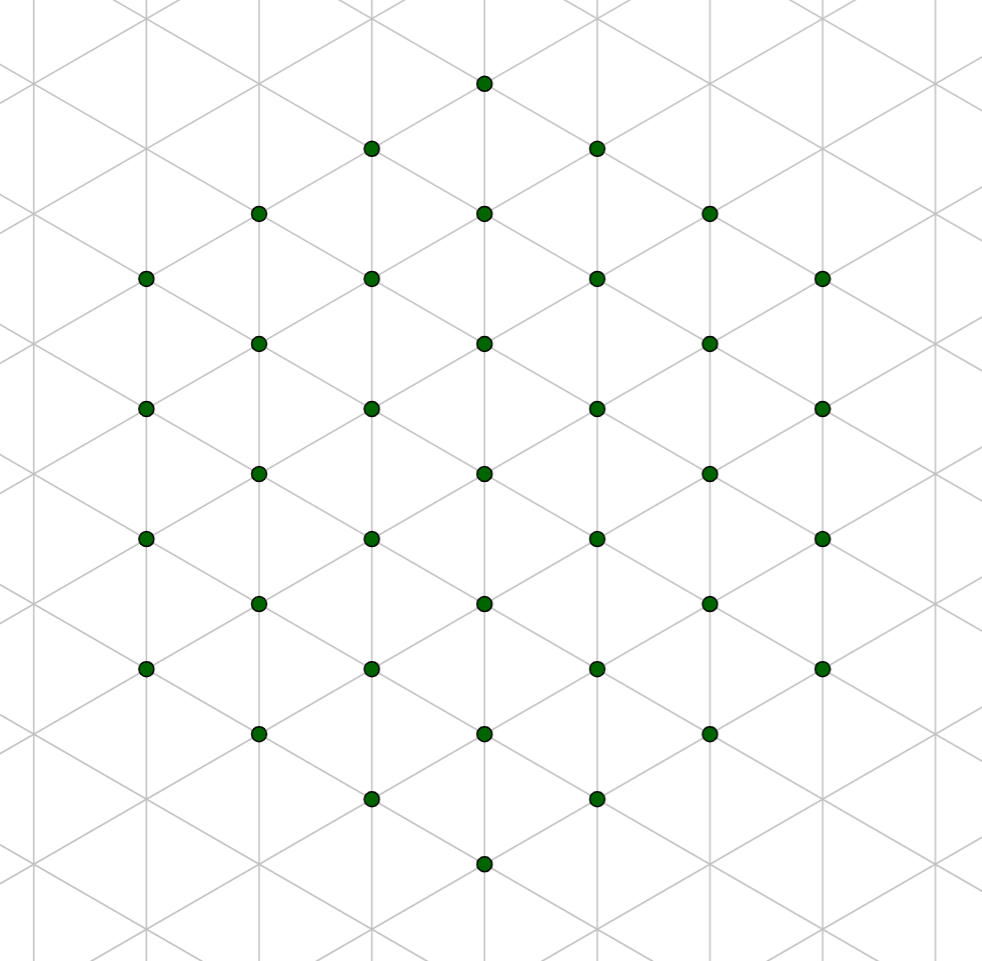}\\
$14$-distance & $15$-distance
\end{tabular}
\end{center}

\begin{rem}
    When Erdős gives the conjecture on maximum $m$-distance sets being on triangular coordinates, he gives the construction of $7$ to $13$-distance set in \cite{erdHos1996maximum}. However, the $10$-distance set construction he provided (with $25$ points) actually contains $11$ distances.
\end{rem}

\begin{center}
\begin{tabular}{cc}
\includegraphics[width=0.3\linewidth]{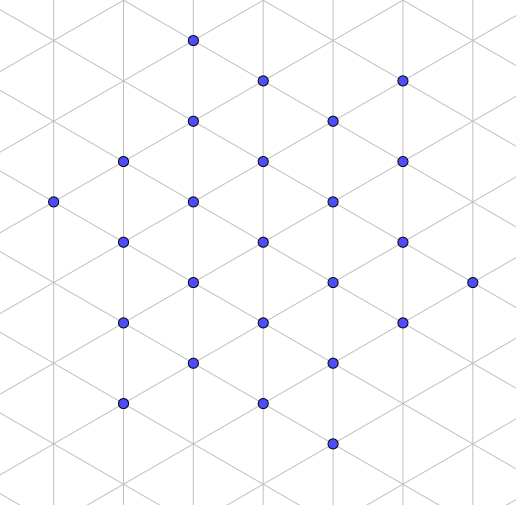}&
\includegraphics[width=0.3\linewidth]{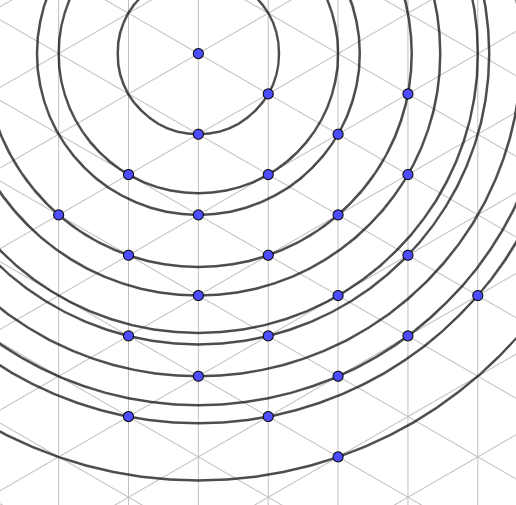}\\
$10$-distance by Erdős & $11$ distances in construction
\end{tabular}
\end{center}

\subsection{$m=17$ to $24$}
Here we present the figures of the construction from $m=17$ to $24$. Here, the code is consistently finding constructions of $m+1$-distance sets distinct from $m$-distance sets, even if they have the same number of points.\newline

\begin{cons}[$17$-distance, $42$ points]
    The figure is symmetric but with only one axis. It did not contain the entire $m=15$ or $m=16$ construction as well. This is also an equiangular hexagon with side lengths alternating $3$ and $4$ with $6$ corners all removed.
\end{cons}
\begin{cons}[$18$-distance, $45$ points]
    This is done by adding eight points to the $m=15$ construction, where a missing point leaves it asymmetrical. This is also an equiangular hexagon with side lengths alternating $3$ and $4$ with $3$ adjacent corners removed.
\end{cons}
\begin{cons}[$19$-distance, $45$ points]
    This is done by adding eight points to the $m=15$ construction as well. However, this addition of points is distributed on four sides instead of three, so the figure is symmetrical. Note that the whole equiangular hexagon does not appear here, so it does not contain the construction of $m=18$ inside.
\end{cons}

\begin{center}
\begin{tabular}{ccc}
\includegraphics[width=0.3\linewidth]{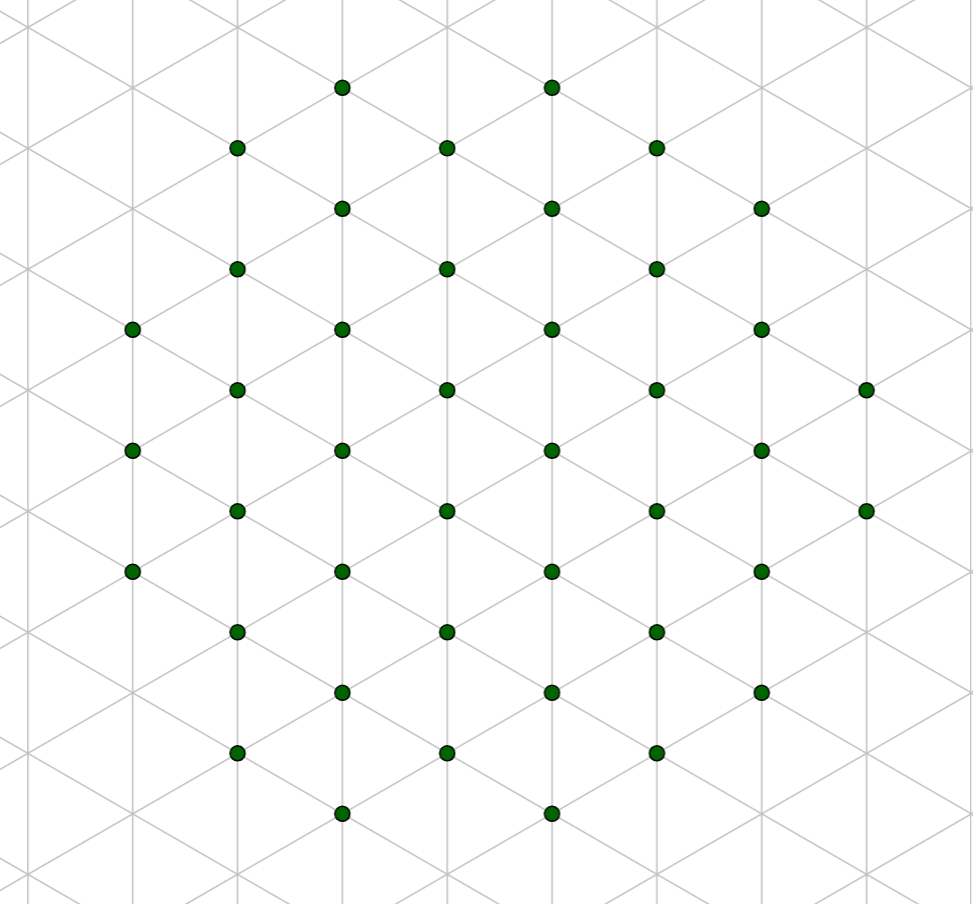}&
\includegraphics[width=0.3\linewidth]{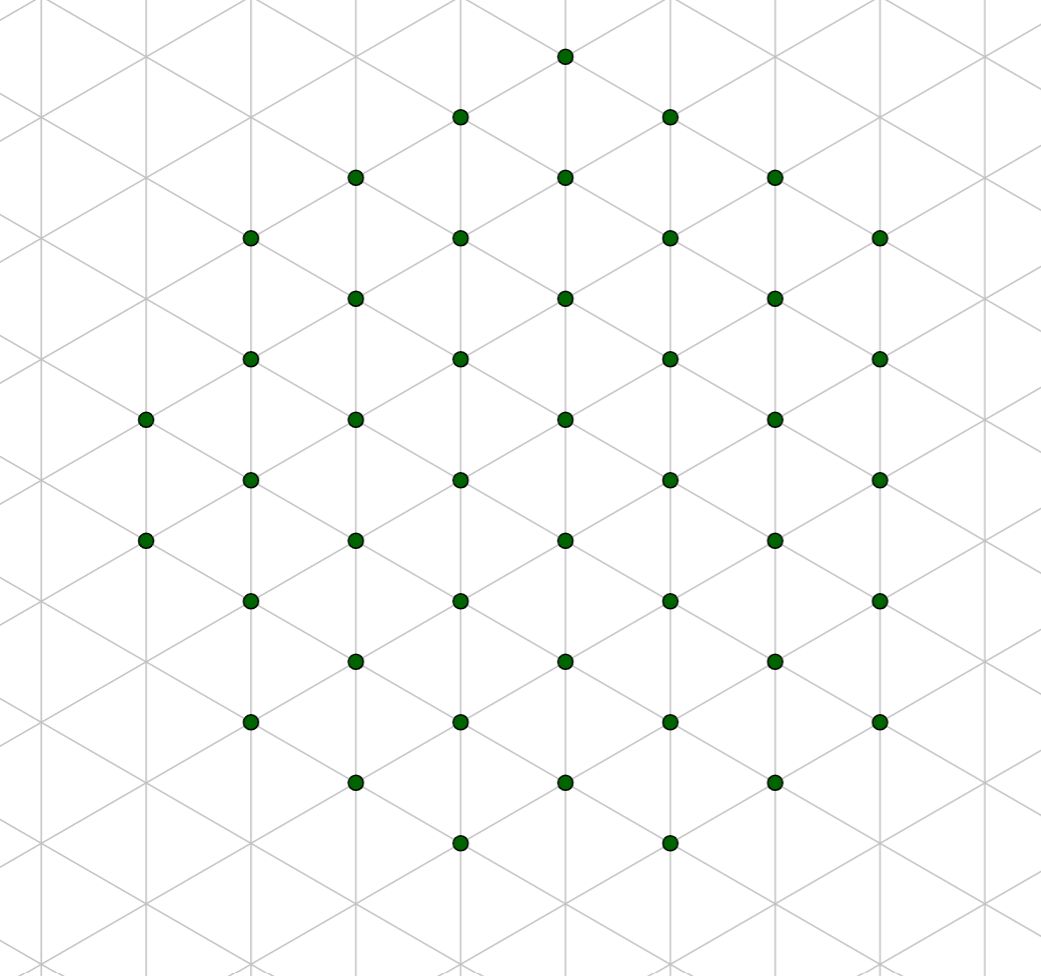}&
\includegraphics[width=0.3\linewidth]{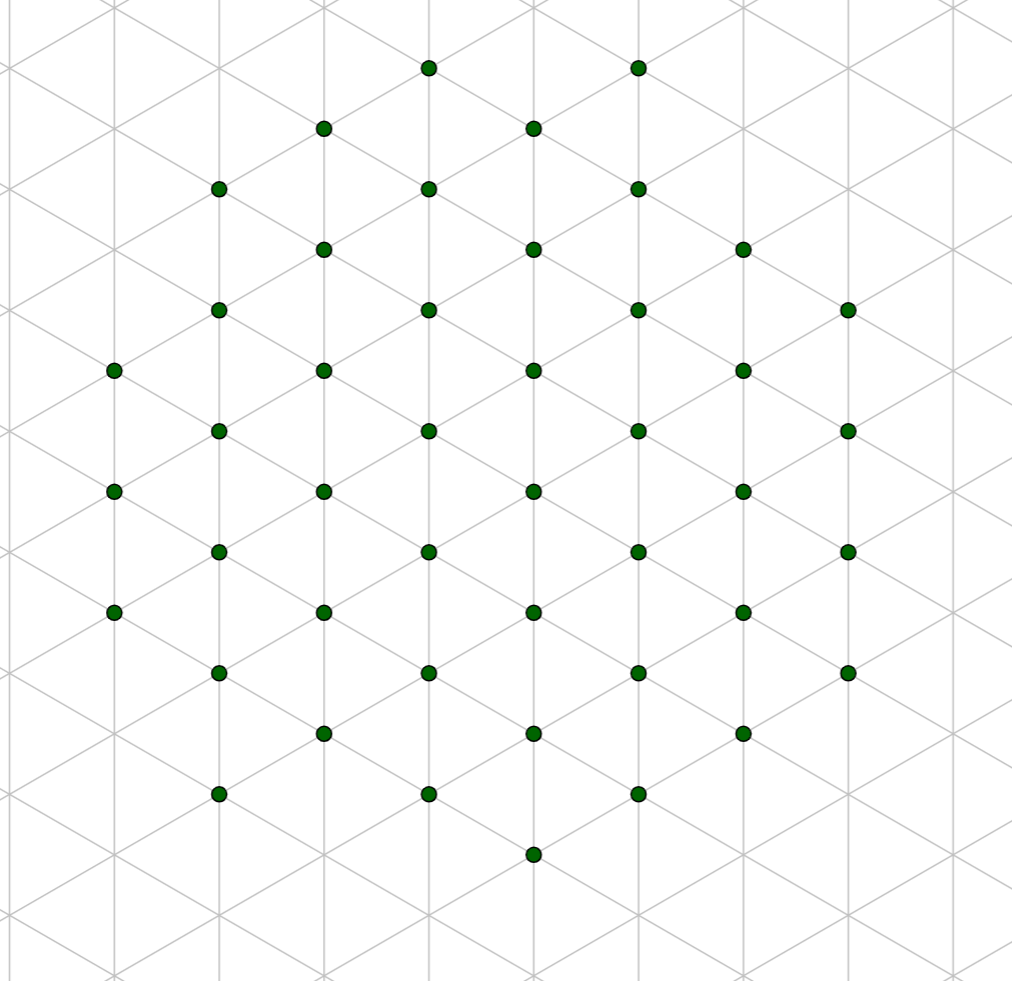}\\
$17$-distance & $18$-distance & $19$-distance
\end{tabular}
\end{center}

\begin{cons}[$20$-distance, $49$ points]
    This is another asymmetrical construction. Removing $12$ points from a regular hexagon with side length $4$.
\end{cons}
\begin{cons}[$21$-distance, $55$ points]
    This one is symmetrical, as it is just removing $6$ corners from a regular hexagon with side length $4$. Giving it $6$ axes of symmetry.
\end{cons}
\begin{cons}[$22$-distance, $58$ points]
    This construction is adding $3$ non-adjacent corner to the $m=21$ construction, or removing $3$ corners from a regular hexagon with side length $4$. Therefore, it has $3$ axes of symmetry.
\end{cons}

\begin{center}
\begin{tabular}{ccc}
\includegraphics[width=0.3\linewidth]{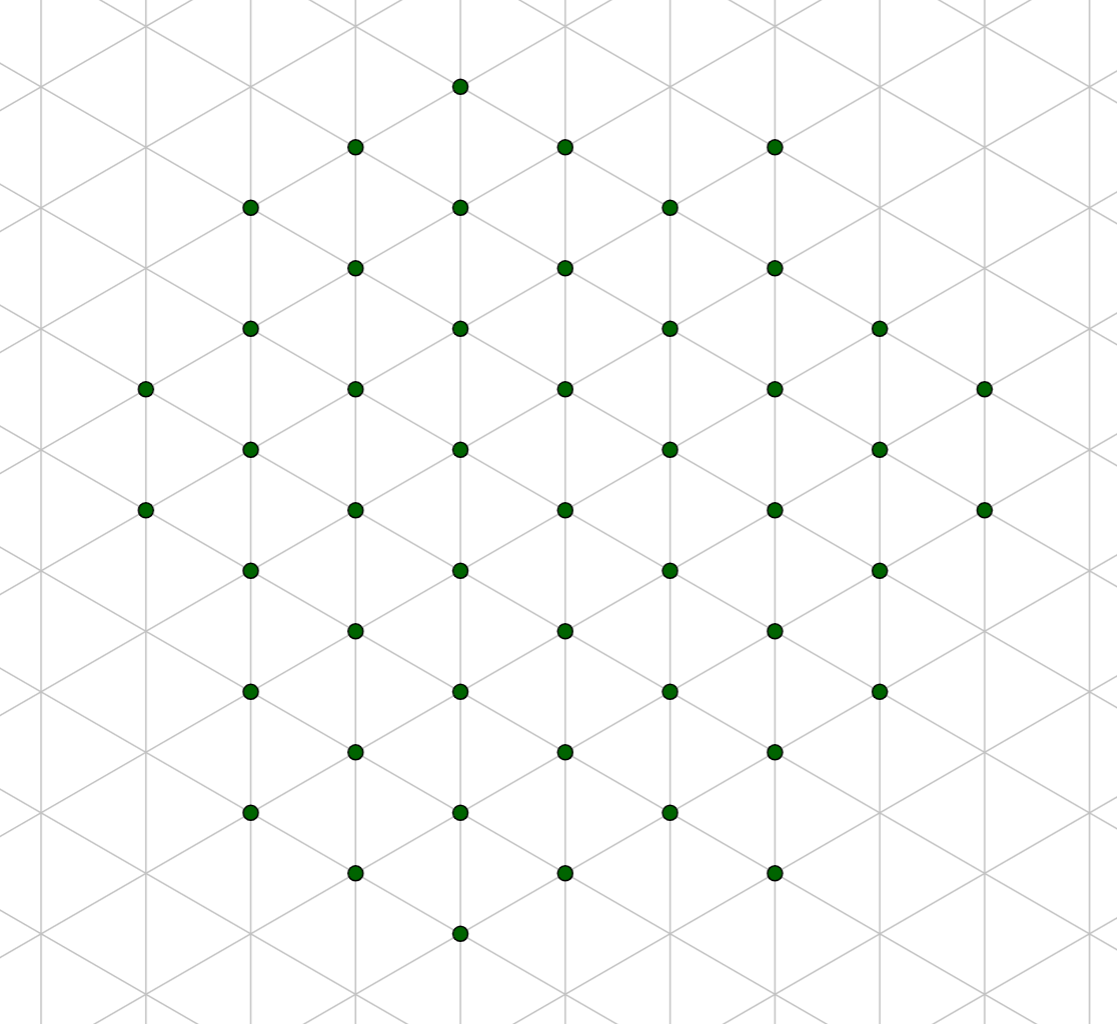}&
\includegraphics[width=0.3\linewidth]{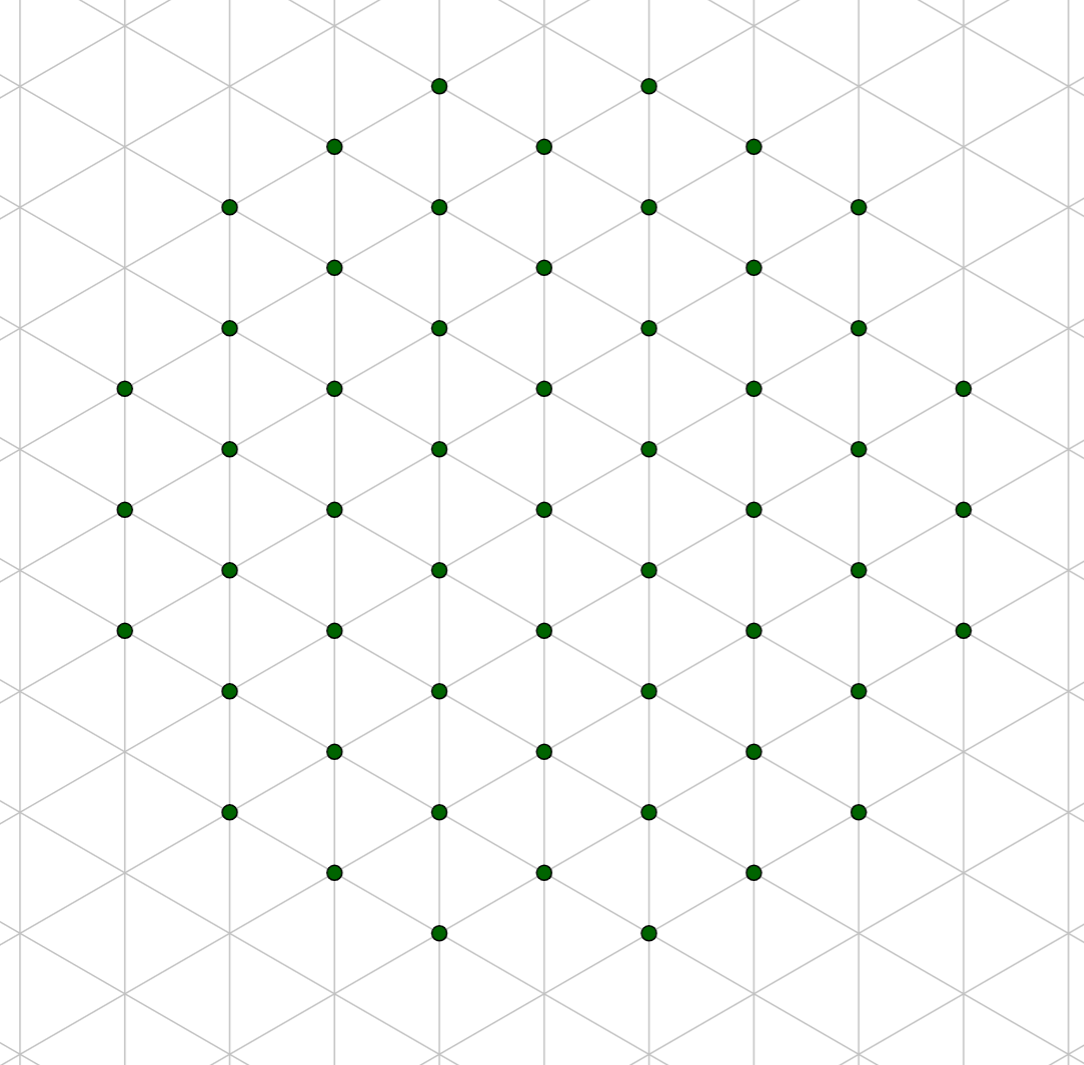}&
\includegraphics[width=0.3\linewidth]{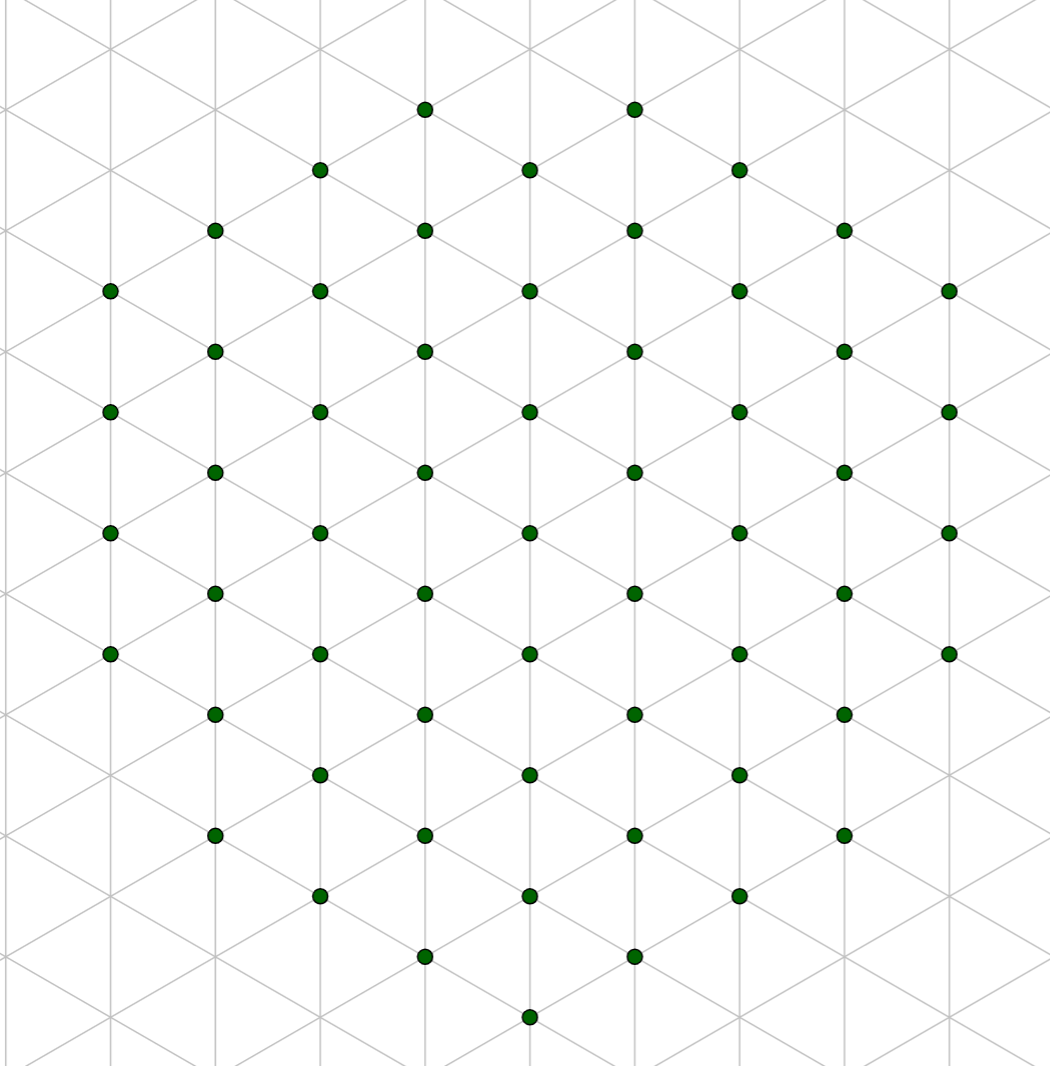}\\
$20$-distance & $21$-distance & $22$-distance
\end{tabular}
\end{center}

\begin{cons}[$23$-distance, $58$ points]
    The code constructed a different figure. This one, however, is asymmetrical. Also, the whole regular hexagon does not show up here.
\end{cons}
\begin{cons}[$24$-distance, $63$ points]
    The construction is still asymmetrical. However, it exceeded $61$ points, which means a maximum $m$-distance set is never a hexagon with side length $5$ when considering the $m$ smallest distances.
\end{cons}

\begin{center}
\begin{tabular}{cc}
\includegraphics[width=0.3\linewidth]{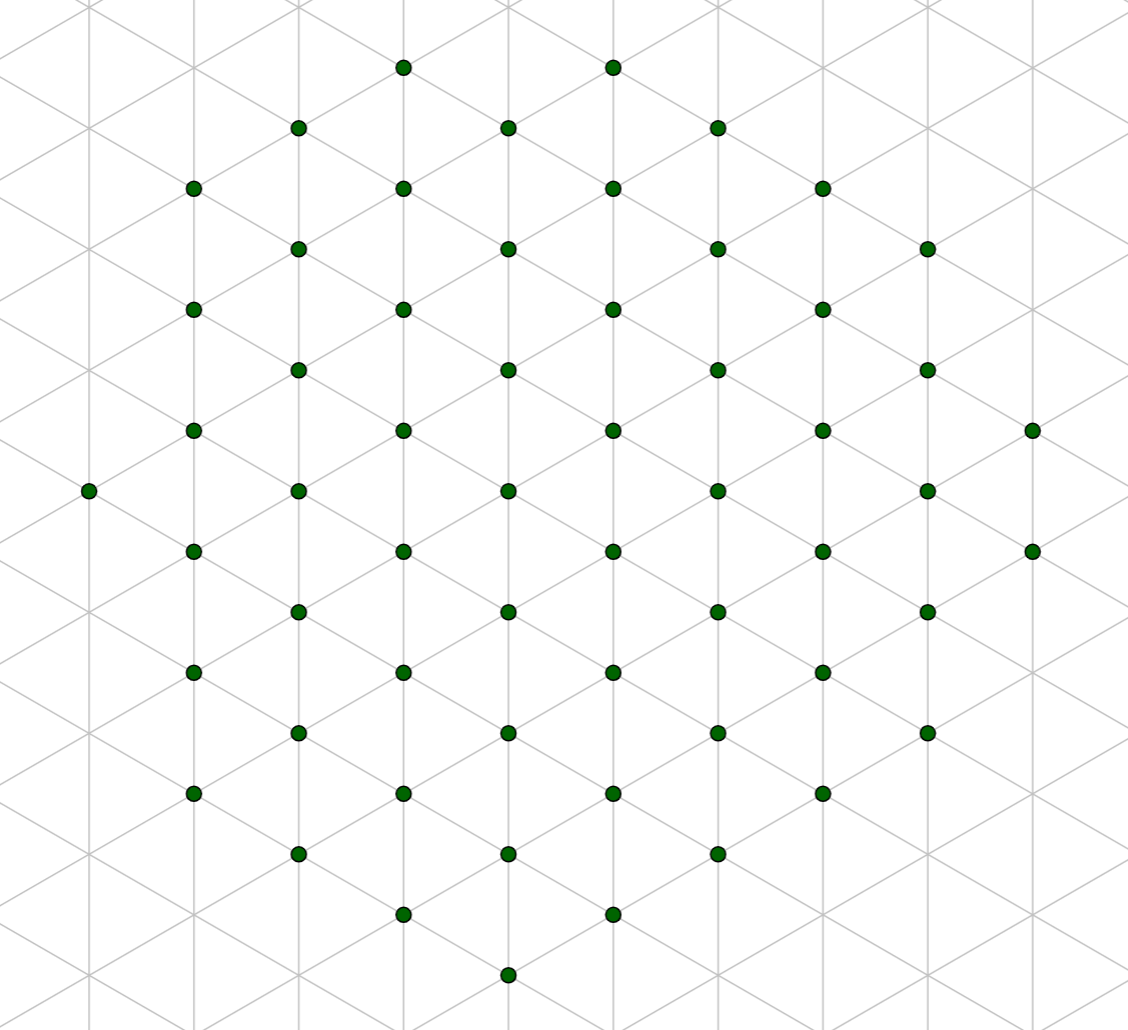}&
\includegraphics[width=0.3\linewidth]{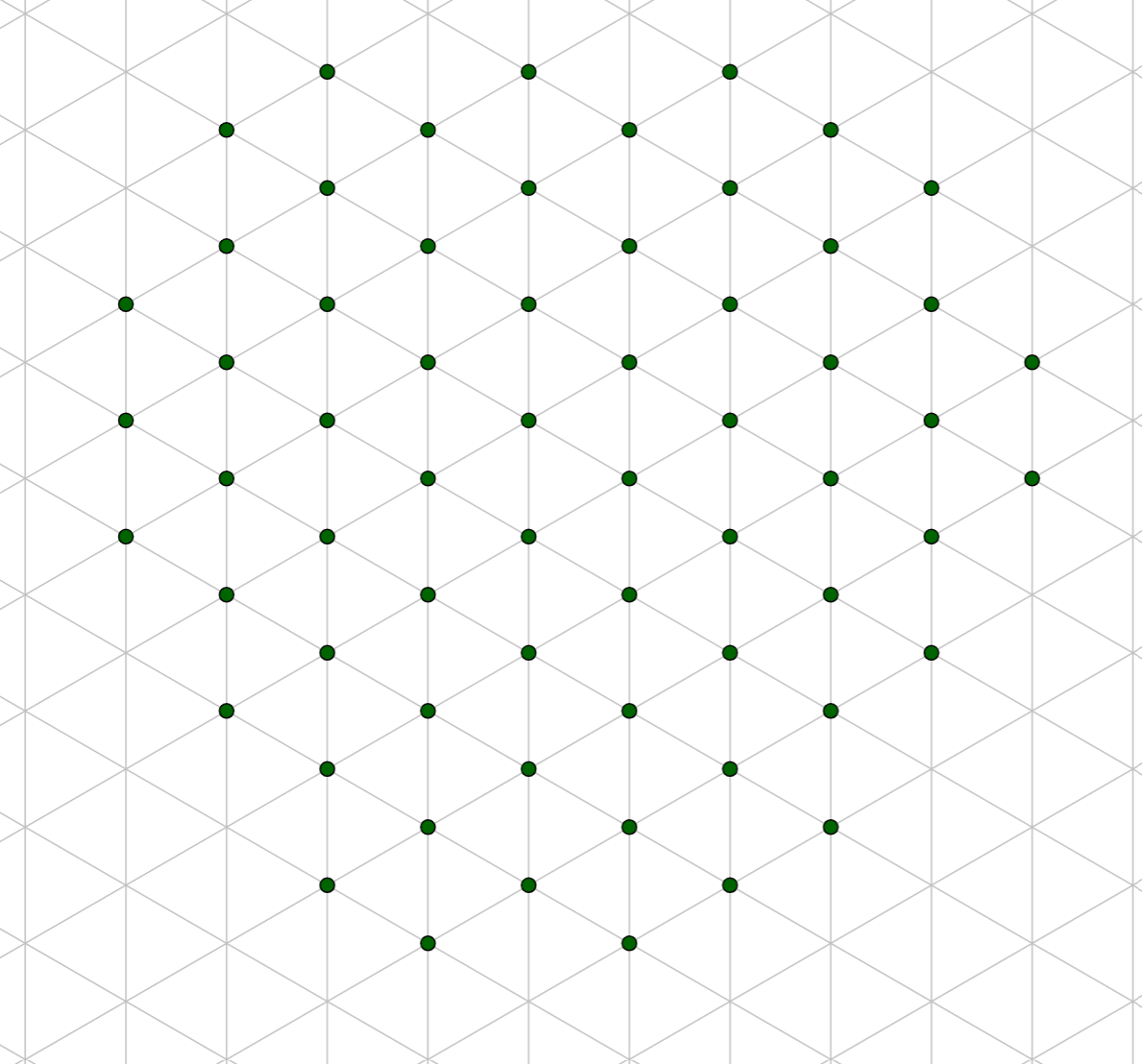}\\
$23$-distance & $24$-distance
\end{tabular}
\end{center}

\subsection{Table for Constructions With the Smallest Distances}
These are the values we got by running the algorithm, with the distances taken to be the $m$ smallest distances on triangular lattice points.
\begin{center}
\begin{tabular}{|c|c|c|}
\hline
    $m$ & size & remark\\ \hline
    $7$ & $16$* & delete three corners from $m=8$\\ \hline
    $8$ & $19$ & regular hexagon side length $2$\\ \hline
    $9$ & $21$ & delete $6$ corners from $m=11$\\ \hline
    $10$ & $24$* & delete $3$ corners from $m=11$\\ \hline
    $11$ & $27$ & equiangular hexagon, side alternating $2, 3$\\ \hline
    $12$ & $27$ & same as $m=11$\\ \hline
    $13$ & $31$ & deleting $6$ corners from regular hexagon\\ \hline
    $14$ & $34$* & deleting $3$ corners from regular hexagon\\ \hline
    $15$ & $37$ & regular hexagon side length $3$\\ \hline
    $16$ & $37$ & same as $m=15$\\ \hline
    $17$ & $42$ & deleting $6$ corners from equiangular hexagon \\\hline
    $18$ & $45$ & deleting $3$ corners from equiangular hexagon\\ \hline
    $19$ & $45$** & symmetric, different from $m=18$\\ \hline
    $20$ & $49$ & delete $12$ points from regular hexagon\\ \hline
    $21$ & $55$ & delete $6$ corners from regular hexagon\\ \hline
    $22$ & $58$* & delete $3$ corners from regular hexagon\\ \hline
    $23$ & $58$** & not symmetric, different from $m=22$\\ \hline
    $24$ & $63$ & not symmetric \\ \hline
\end{tabular}
\end{center}
Sizes with single star mean that confirmed alternate constructions exist with the smallest $m$-distances. Sizes with double stars mean other distances (not the $m$ smallest ones) would give different constructions, and hence the construction here with double stars are not maximum distance sets.

\subsection{Discussions}
As we can see, many of the construction are based on an equiangular hexagon, and removing some of the corners from it. Some observations on the removals give us some interesting results.
\begin{pp}
    The algorithm removed $3$ corners for $7$-distance, $10$-distance, $14$-distance, $18$-distance, and $22$-distance sets. However, the removals are a bit different. The $7$-distance, $10$-distance, and $22$-distance have none of the three removed corner being adjacent. Which is not the case for $14$-distance and $18$-distance. However, this structure can also be applied to $7, 10, 22$-distance as well.
\end{pp}
\begin{center}
\begin{tabular}{ccc}
\includegraphics[width=0.3\linewidth]{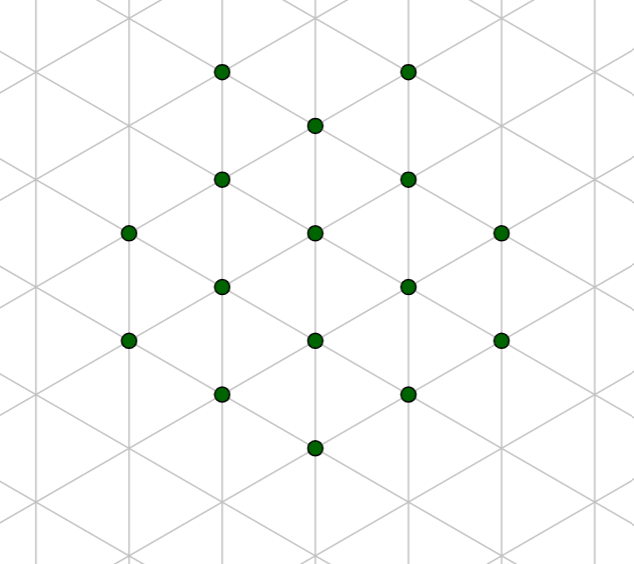}&
\includegraphics[width=0.3\linewidth]{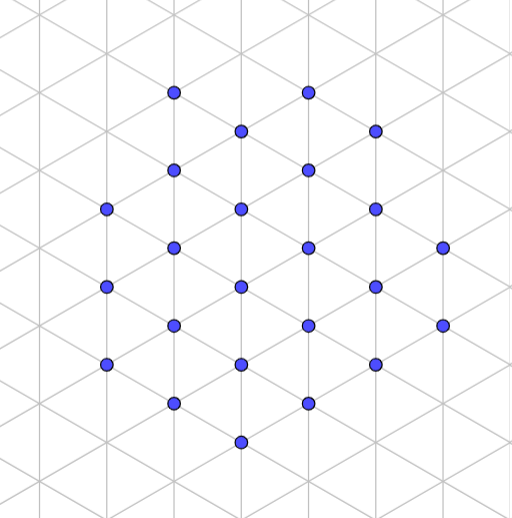}&
\includegraphics[width=0.3\linewidth]{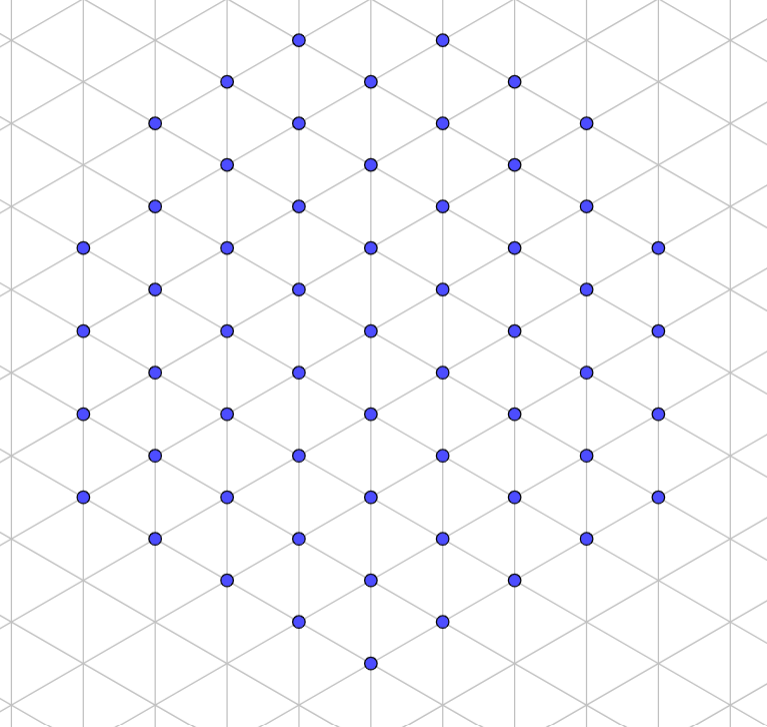}\\
$7$-distance & $10$-distance & $22$-distance
\end{tabular}
\end{center}
\begin{pp}
    Similar to the last property. $14$-distance and $18$-distance can also be constructed with the $7$-distance-like structure.
\end{pp}
\begin{pp}
    The $20$-distance set the construction gave is deleting $12$ points. However, the deleted points look a bit arbitrary and not symmetric. There is another construction that is more "structured" and symmetric in all six corners.
\end{pp}
\begin{center}
\begin{tabular}{ccc}
\includegraphics[width=0.3\linewidth]{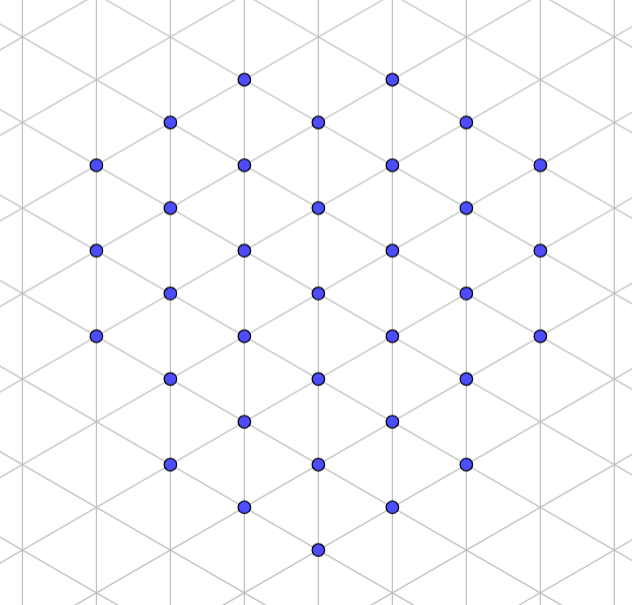}&
\includegraphics[width=0.3\linewidth]{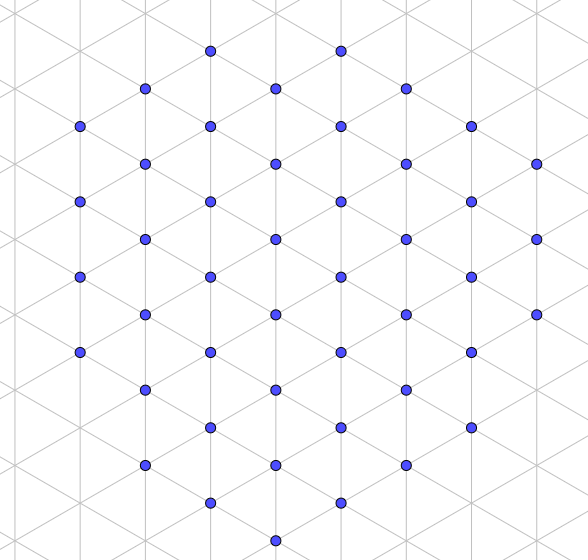}&
\includegraphics[width=0.3\linewidth]{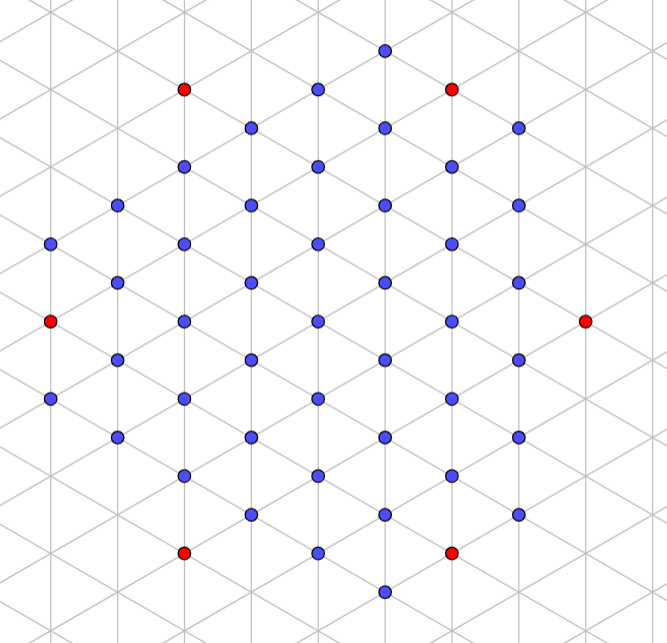}\\
$14$-distance & $18$-distance & $20$-distance
\end{tabular}
\end{center}

Finally, there are the equiangular hexagons that do not remove any of the corners. Those will be discussed in the next section

\section{Construction With Other Distances}
\subsection{Some Observations on Hexagons}
We noticed that $15$-distance featured a regular hexagon, and $14$-distance has $3$ corner removed from it. However, the same structure is not seen with $18$-distance and $19$-distance, or with $22$-distance and $23$-distance. Why? Because these hexagons do not contain some of the smallest $19$ (or $23$) distances on triangular lattices.
\begin{pp}
    An equiangular hexagon with side length alternating $3$ and $4$ is a $19$-distance set. However, it does not contain the distance $\sqrt{48}$, which is the $19$th smallest triangular coordinate distance. Instead, it contains $7$, which is the $20$th smallest. This construction has $48$ points, which is more than the $45$ given in the previous section.
\end{pp}
\begin{pp}
    A regular hexagon with side length $4$ is a $23$-distance set. However, it does not contain the distance $\sqrt{61}$, which is the $23$rd smallest triangular coordinate distance. Instead, it contains $8$, which is the $25$th smallest. This construction has $63$ points, which is more than the $58$ given in the previous section.
\end{pp}
\begin{center}
\begin{tabular}{cc}
\includegraphics[width=0.3\linewidth]{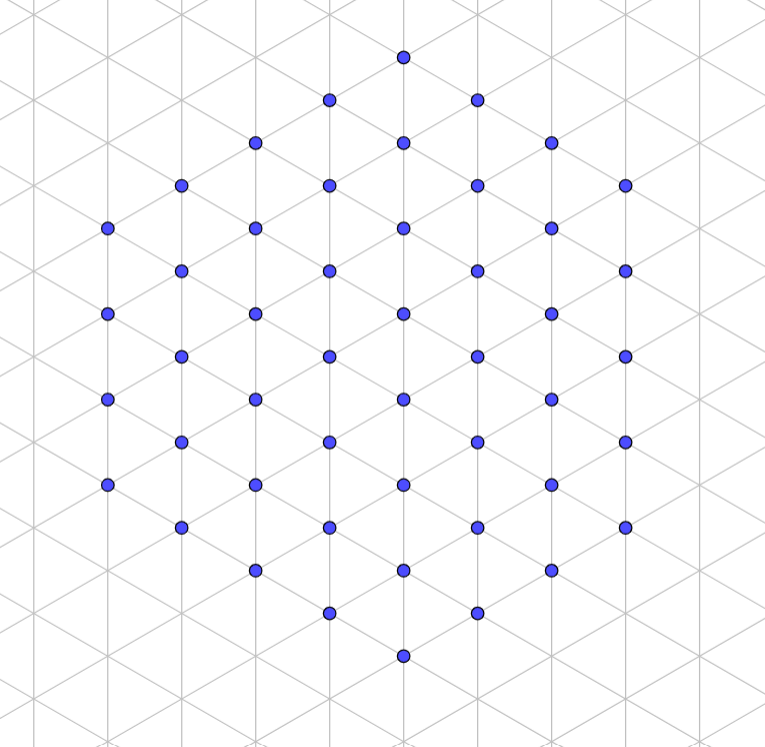}&
\includegraphics[width=0.3\linewidth]{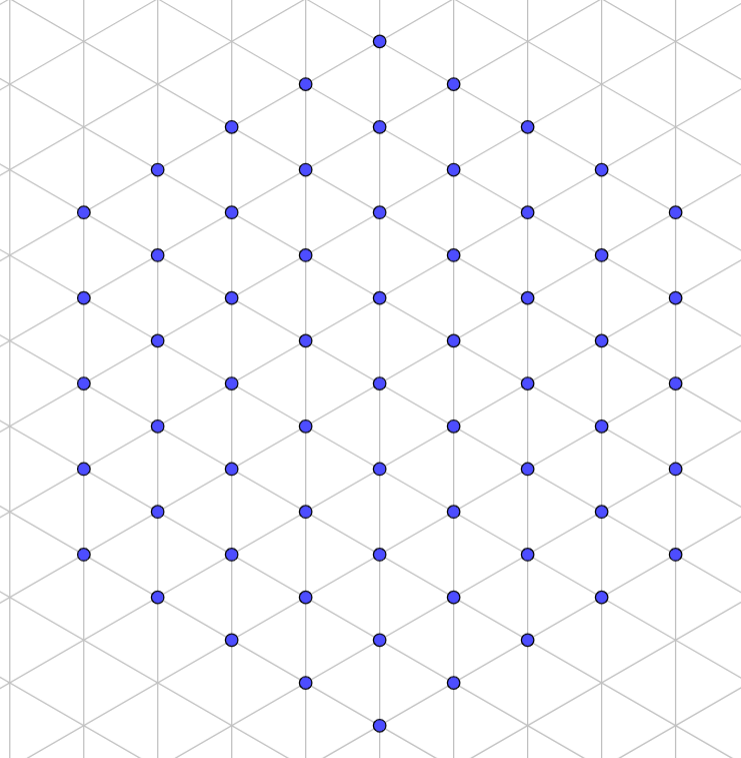}\\
$19$-distance & $23$-distance
\end{tabular}
\end{center}
For convenience, we make the following definition, and will use it in the next subsections.
\begin{df}
    We call an equiangular hexagon with alternating side lengths $k$ and $k+1$ a $(k, k+1)$-equiangular hexagon. 
\end{df}

\subsection{Multiplicity of Each Distance in Constructions}
When running the maximum clique algorithm, each of the $m$ distances in $\mathcal S$ would count as an edge. However, some might appear more often. In this section, we will see how often does each distance appear in such constructions.\newline
\begin{df}[Multiplicity Array]
    A multiplicity array of a construction with at most $m$ distances is an array of length $m$. The $i$th entry is the multiplicity of the $i$th smallest distance in a construction.
\end{df}
First, we will see the multiplicity array for constructions of $m\leq 16$. We will use squares of distances to represent the distances on the table.

\begin{center}
\begin{tabular}{|c|c|c|c|c|c|c|c|c|c|c|c|c|c|c|c|c|}
\hline
    $m$ & $1$ & $3$ & $4$ & $7$ & $9$ & $12$ & $13$ & $16$ & $19$ & $21$ & $25$ & $27$ & $28$ & $31$ & $36$ & $37$ \\\hline
    $7$ & $33$ & $24$ & $18$ & $24$ & $9$ & $6$ & $6$ &&&&&&&&&\\\hline
    $8$ & $42$ & $30$ & $27$ & $36$ & $12$ & $9$ & $12$ & $3$ &&&&&&&&\\\hline
    $9$ & $45$ & $36$ & $30$ & $42$ & $18$ & $9$ & $18$ & $6$ & $6$ &&&&&&&\\\hline
    $10$ & $54$ & $42$ & $36$ & $54$ & $21$ & $15$ & $27$ & $9$ & $12$ & $6$ &&&&&& \\\hline
    $11$ & $63$ & $48$ & $45$ & $66$ & $27$ & $21$ & $36$ & $12$ & $18$ & $12$ & $3$ &&&&& \\\hline
    $12$ & $63$ & $48$ & $45$ & $66$ & $27$ & $21$ & $36$ & $12$ & $18$ & $12$ & $3$ & $\mathbf{0}$ &&&&\\\hline
    $13$ & $72$ & $60$ & $51$ & $84$ & $36$ & $27$ & $48$ & $21$ & $24$ & $24$ & $6$ & $6$ & $6$ &&&\\\hline
    $14$ & $81$ & $66$ & $60$ & $96$ & $41$ & $33$ & $60$ & $24$ & $36$ & $30$ & $9$ & $7$ & $12$ & $6$ && \\\hline
    $15$ & $90$ & $72$ & $69$ & $108$ & $48$ & $39$ & $72$ & $27$ & $48$ & $36$ & $12$ & $12$ & $18$ & $12$ & $3$ & \\\hline
    $16$ & $90$ & $72$ & $69$ & $108$ & $48$ & $39$ & $72$ & $27$ & $48$ & $36$ & $12$ & $12$ & $18$ & $12$ & $3$ & $\mathbf{0}$ \\\hline
\end{tabular}
\end{center}

Next, we would analyze some cases where the hexagonal construction yields more points than the max clique. Starting with $m=19$, the smallest case where the hexagon construction overtakes the maximum clique algorithm.

\begin{pp}
    The $m=19$ construction using maximum clique algorithm has multiplicity array $\{110, 93, 87, 144, 64, 54, 102, 43, 72, 62, 24, 21, 40, 30, 9, 16, 12, 6, 1\}$, while the construction with the hexagon has multiplicity array $\{120, 99, 96, 156, 72, 60,\\ 114, 48, 84, 72, 27, 27, 48, 36, 12, 24, 18, 12, 3\}$. Where the $19$th distance in the second construction is $7$ instead of $\sqrt{48}$
\end{pp}

We found that the maximum clique construction has a distance ($\sqrt{48}$) with multiplicity only $1$. Next up is $m=23$, where the hexagon is also bigger by $3$ points.

\begin{pp}
    The $m=23$ construction using maximum clique algorithm has multiplicity array $\{146, 126, 119, 204, 93, 81, 154, 68, 118, 109, 45, 40, 78, 67, 26, 47, 41, 32,\\ 10, 27, 14, 6, 2\}$, while the construction with the hexagon has multiplicity array $\{156, 132,\\ 129, 216, 102, 87, 168, 75, 132, 120, 48, 48, 90, 72, 27, 60, 48, 36, 15, 36, 18, 12, 3\}$. Where the $23$rd distance in the second construction is $8$ instead of $\sqrt{61}$.
\end{pp}

We found that $\sqrt{61}$ has multiplicity only $1$. Finally, as $m=28$ also has a hexagon construction, we will also analyze its multiplicity array for the two methods.
\newpage
\begin{pp}
    The $m=28$ construction using maximum clique algorithm has multiplicity array $\{186, 162, 156, 270, 126, 111, 216, 97, 174, 162, 70, 66, 126, 111, 45, 90,\\ 81, 66, 27, 75, 42, 30, 22, 18, 9, 12, 6, 0\}$, while the construction with the hexagon has multiplicity array $\{195, 168, 165, 282, 135, 117, 228, 105, 186, 174, 75, 72, 138, 120, 48, 102,\\ 90, 72, 33, 87, 48, 36, 30, 24, 12, 18, 12, 3\}$. Where the $23$rd distance in the second construction is $9$ instead of $\sqrt{75}$.
\end{pp}
These examples give an explanation on why sometimes the $m$ smallest distances would not perform as well as the hexagon construction. Some larger distances in such constructions are often "underused" to the point that the $m$ smallest distances and the $m-1$ smallest distances yield the same number of points in their constructions. The largest distance in those constructions usually contains a distance with multiplicity $2$ or less. These hexagon constructions take out the largest distance (by deleting only $1$ or $2$, or even no points), and add more points for a larger construction.
Therefore, what hexagons do is "kick" the distances with low multiplicity, and replace them with the nearest bigger integer distance. This gives a conjecture about maximum $m$-distance sets would never have a distance with multiplicity $2$ or lower.
Also, there is an additional observation that can be made here.
\begin{pp}
    If the construction is based on the hexagon construction, then every entry of the multiplicity array would be a multiple of $3$.
\end{pp}
This can be easily seen by the symmetry of the hexagon and the removed points.



\subsection{Constructions with Hexagons}
First, we will see why we run the previous algorithm with the smallest $19$ distances, and can not construct the hexagons for $m=19$.

\begin{pp}
    $(k,k+1)$-equiangular hexagons and regular hexagons cannot be constructed with the smallest $m$-distances when $m\geq 19$.
\end{pp}
\begin{proof}
    When $m=19$, the longest diagonal has length $7$, so we would consider the case that the length $L$ of the diagonal is at least $7$. The largest distance in this hexagon is $L$, and the second largest is $\sqrt{(L-\frac12)^2+\frac34}=\sqrt{L^2-L+1}$. Note that the length $\sqrt{(L-\frac12)^2+(\frac{3\sqrt3}{2})^2}=\sqrt{L^2-L+8}$ is less than $L$ when $L\geq 8$. Therefore, for $m\geq 23$, the hexagon would miss a distance between its two largest distances. $m=19$ can be checked by the picture in the previous subsections.
\end{proof}

Therefore, we might consider to construct $m$-distance sets by the shape of hexagons.
\begin{pp}[Hexagon Constructions]
    If we want to build an $m$-distance set, we find the smallest $m'\geq m$ such that there exists a $(k, k+1)$-equiangular hexagon or regular hexagon forming an $m'$-distance set, and repeatedly remove one of the points on the diameters until it becomes an $m$-distance set. We would also consider the smaller hexagon to check if it has more points than the previous construction done by removing points. For example, if we construct a $29$-distance set by removing points from the regular hexagon of $34$-distance set, then we would have fewer points than the equiangular hexagon of $28$-distance set. In this case, we should take the smaller hexagon instead of removing points from the larger hexagon.
\end{pp}
We will compare the maximum clique method with the hexagon construction method. To make the maximum clique algorithm work faster, we do the following: assuming that $(m-1)$-distance set gets a maximum clique of $t$ vertices, as a vertex in a clique of size $t$ has degree at least $t-1$, we remove all vertices with degree less than $t-1$ when we construct $m$-distance set.
Here are the results of comparing the two methods. The value shows the maximum constructions with at most $m$ distances.
\begin{center}
\begin{tabular}{|c|c|c|c|}
\hline
    $m$ & max clique& hexagon & remark\\ \hline
    $23$ & $58$ & $\mathbf{61}$ & regular hexagon side length $4$\\ \hline
    $24$ & $63$ & $63$ & remove $12$ points from $4, 5$-equiangular\\ \hline
    $25$ & $63$ & $63$ & same as $m=24$ with the hexagon construction\\ \hline
    $26$ & $69$ & $69$ & delete $6$ corners from $4, 5$-equiangular\\ \hline
    $27$ & $72$ & $72$ & delete $6$ corners from $4, 5$-equiangular\\ \hline
    $28$ & $72$ & $\mathbf{75}$ & $4, 5$-equiangular hexagon\\ \hline
    $29$ & $73$ & $\mathbf{75}$ & same as $m=28$\\ \hline
    $30$ & $79$ & $79$ & deleting $12$ points from regular hexagon side length $5$\\ \hline
    $31$ & $79$ & $79$ & same as $m=30$ with the hexagon construction\\ \hline
    $32$ & $85$ & $85$ & deleting $6$ corners from regular hexagon side length $5$\\ \hline
    $33$ & $88$ & $88$ & deleting $3$ corners from regular hexagon side length $5$ \\\hline
    $34$ & $90$ & $\mathbf{91}$ & regular hexagon side length $5$\\ \hline
\end{tabular}
\end{center}
\begin{center}
\begin{tabular}{ccc}
\includegraphics[width=0.3\linewidth]{m-distance/2D-23-dis-alt1.png}&
\includegraphics[width=0.3\linewidth]{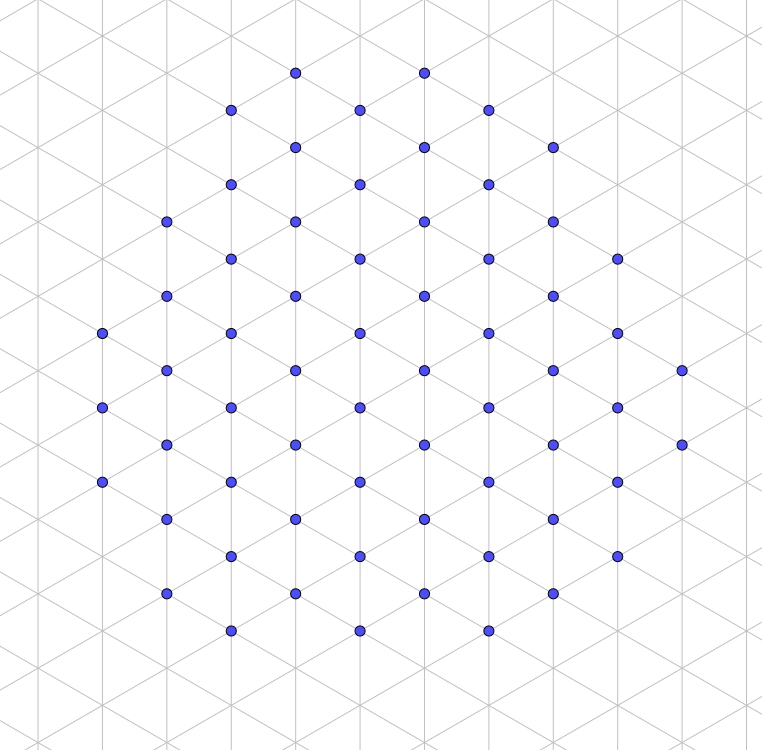}&
\includegraphics[width=0.3\linewidth]{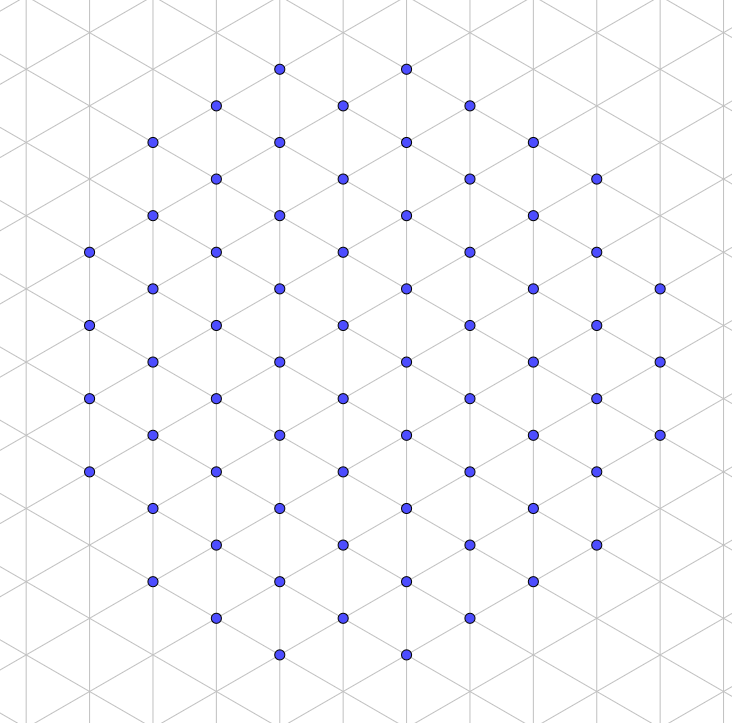}\\
$23$-distance & $24$-distance & $26$-distance
\end{tabular}
\end{center}
\begin{center}
\begin{tabular}{ccc}
\includegraphics[width=0.3\linewidth]{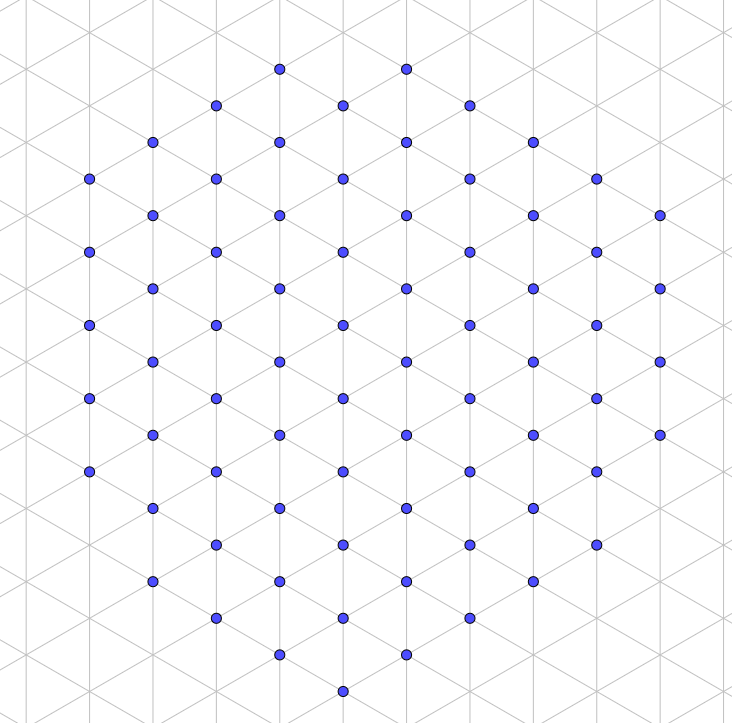}&
\includegraphics[width=0.3\linewidth]{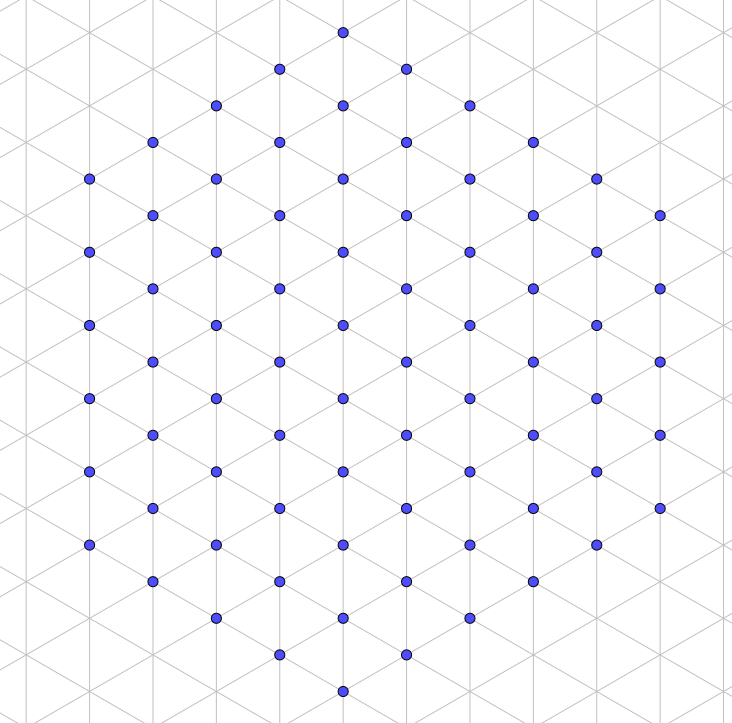}&
\includegraphics[width=0.3\linewidth]{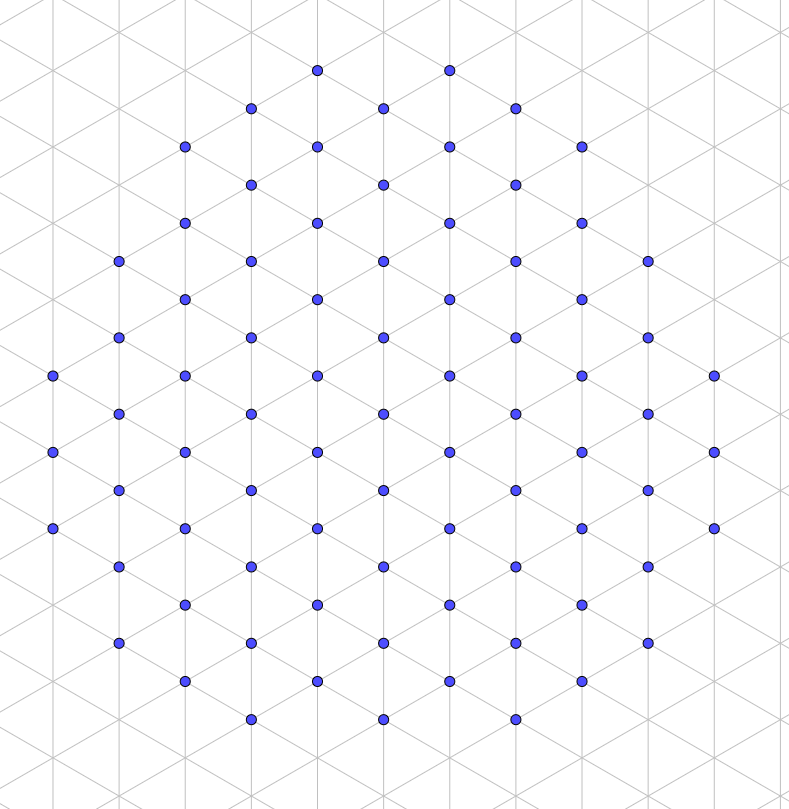}\\
$27$-distance & $28$-distance & $30$-distance
\end{tabular}
\end{center}
\begin{center}
\begin{tabular}{ccc}
\includegraphics[width=0.3\linewidth]{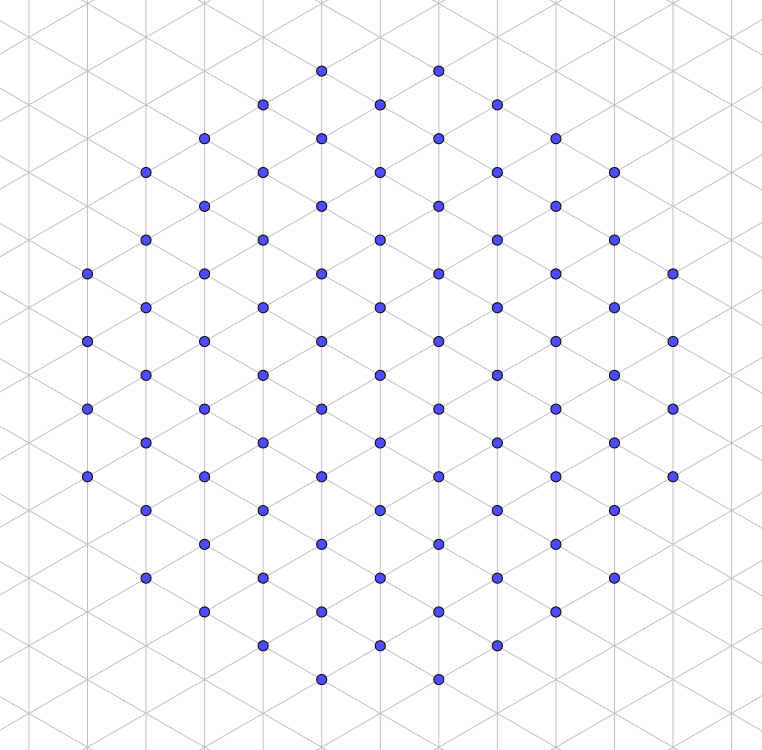}&
\includegraphics[width=0.3\linewidth]{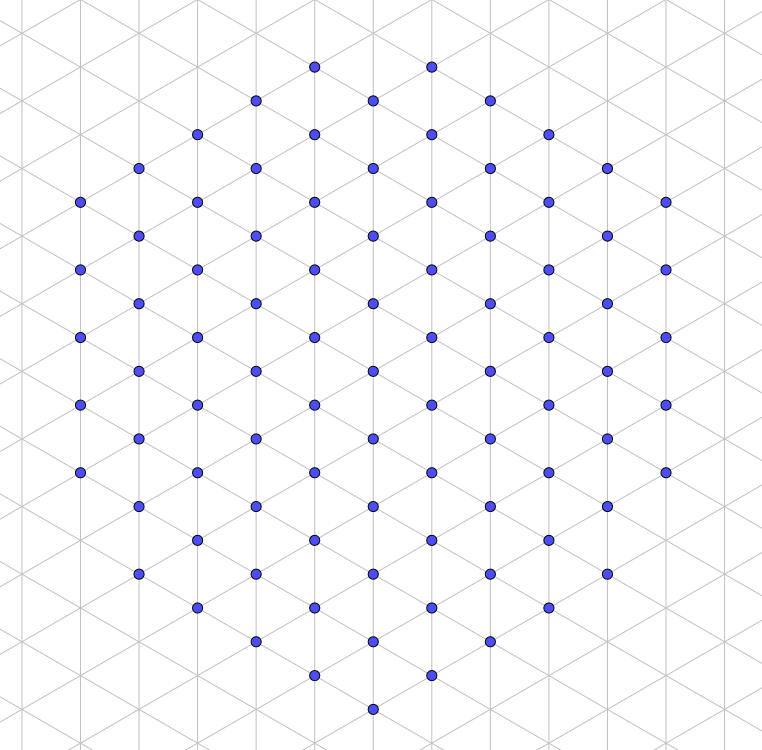}&
\includegraphics[width=0.3\linewidth]{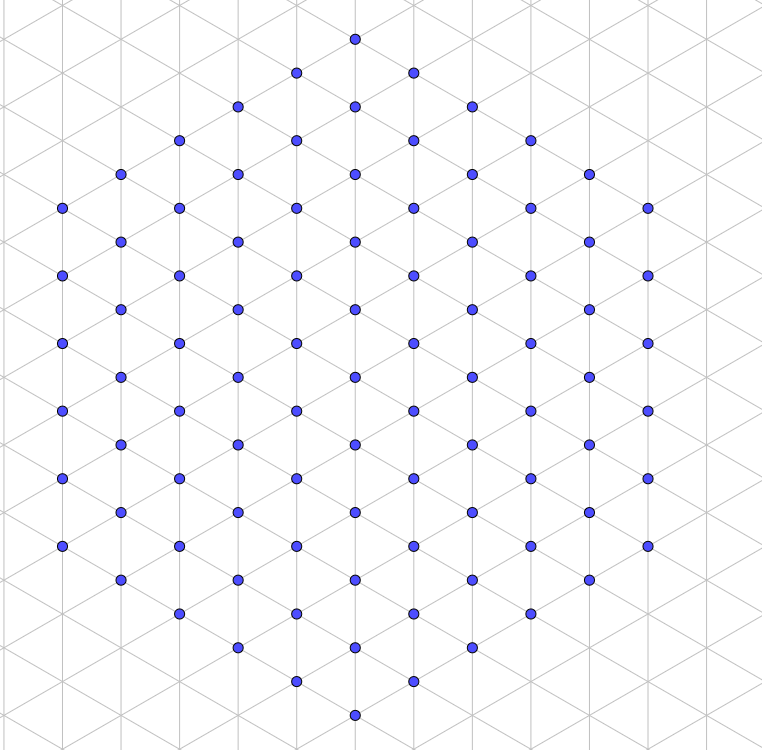}\\
$32$-distance & $33$-distance & $34$-distance
\end{tabular}
\end{center}

Therefore, we conjecture that the construction in the above table of hexagon would be the maximum.

\section{Conclusion and Further Discussions}
\subsection{Conclusion on the Two Methods}
\begin{con}
    The construction by the method of maximum clique with the first $m$ smallest distances gives the largest construction, except for $m=19, 23, 28, 29, 34$. Therefore, the clique method works for most of the cases.  
\end{con}
\begin{con}
    The construction by using the first $m$ smallest distances in triangular lattices can not attain the maximum size for $m=19, 23, 28, 29, 34$,  since the construction of hexagons can have more points. 
\end{con}
These ones often have the maximum clique algorithm finding the same answer with the $m-1$ smallest distances, because the multiplicity of the $m$th distances are all $3$ or less.

When the $m$-distance set construction of hexagon is larger than the construction by maximum clique, the construction of maximum clique for $m$-distance set often has the same size as $(m-1)$-distance set since the largest distance has small multiplicity. 

When we compare our constructions to the ones in \cite{balaji2019lattice}, we noticed that some of our constructions are larger. For the $18$-distance set, the paper removed points from a regular hexagon with side length $4$ to get a $43$ points construction. On the contrary, we removed three corners from a $(3, 4)$-equiangular hexagon to get a $45$ points construction. There was also a $70$ points construction for $29$-distance set in that paper. However, we were able to construct a set with $75$ points and only $28$ distances by taking a $(4, 5)$-equiangular hexagon. The main reason that we are finding larger construction is mainly about the observations made on $(k, k+1)$-equiangular hexagons, as opposed to only using regular hexagons and hexagons with line symmetry.

\subsection{Further Discussions}
    Comparing the two methods, the construction of hexagons works better, so we conjecture that the construction of hexagons will attain the maximum i.e. all triangular lattices in a regular hexagon, or a $(k, k+1)$-equiangular hexagon will form a maximum $m$-distance set. In addition, there is a stronger version that can be used on every $m$: if $m\geq 3$ then at least one construction of maximum $m$-distance set can be found by repeatedly removing points from a regular hexagon or a $(k, k+1)$-equiangular hexagon. Also, we made a conjecture based on some non-optimal construction presented by the maximum clique algorithm.
    For all $m\geq 3$, any maximum $m$-distance set does not have a distance with multiplicity less than $3$. Finally, for the $7$-distance and $8$-distance construction, we also conjecture that the maximum $7$-distance set in triangular lattices has $16$ points, and the maximum $8$-distance set in triangular lattices has $19$ points. In addition, they are only realized by the constructions found in subsection 2.1 and 2.4. If Erdős's conjecture is true, then these would be the only maximum $7$ and $8$ distance sets in $\mathbb R^2$.



%% file: main.bbl
\begin{thebibliography}{10}

\bibitem{adachi2017maximal}
Saori Adachi, Rina Hayashi, Hiroshi Nozaki, and Chika Yamamoto.
\newblock Maximal m-distance sets containing the representation of the hamming graph h (n, m).
\newblock {\em Discrete Mathematics}, 340(3):430--442, 2017.

\bibitem{ahmeddistance}
Tanbir Ahmed and David~Jacob Wildstrom.
\newblock On distance sets in the triangular lattice.

\bibitem{balaji2019lattice}
Vajresh Balaji, Olivia Edwards, Anne~Marie Loftin, Solomon Mcharo, Lo~Phillips, Alex Rice, and Bineyam Tsegaye.
\newblock Lattice configurations determining few distances.
\newblock {\em arXiv preprint arXiv:1911.11688}, 2019.

\bibitem{bannai2012maximal}
Eiichi Bannai, Takahiro Sato, and Junichi Shigezumi.
\newblock Maximal m-distance sets containing the representation of the johnson graph j (n, m).
\newblock {\em Discrete Mathematics}, 312(22):3283--3292, 2012.

\bibitem{carraghan1990exact}
Randy Carraghan and Panos~M Pardalos.
\newblock An exact algorithm for the maximum clique problem.
\newblock {\em Operations Research Letters}, 9(6):375--382, 1990.

\bibitem{erdHos1996maximum}
Paul Erd{\H{o}}s and Peter Fishburn.
\newblock Maximum planar sets that determine k distances.
\newblock {\em Discrete Mathematics}, 160(1-3):115--125, 1996.

\bibitem{oeis:A003136}
The OEIS~Foundation Inc.
\newblock Sequence a003136, loeschian numbers: numbers of the form $x^2 + xy + y^2$; norms of vectors in a2 lattice.
\newblock \url{https://oeis.org/A003136}.
\newblock [Online; accessed 16-June-2025].

\bibitem{lisonvek1997new}
Petr Lison{\v{e}}k.
\newblock New maximal two-distance sets.
\newblock {\em Journal of Combinatorial Theory, Series A}, 77(2):318--338, 1997.

\bibitem{shinohara2004classification}
Masashi Shinohara.
\newblock Classification of three-distance sets in two dimensional euclidean space.
\newblock {\em European Journal of Combinatorics}, 25(7):1039--1058, 2004.

\bibitem{shinohara2008uniqueness}
Masashi Shinohara.
\newblock Uniqueness of maximum planar five-distance sets.
\newblock {\em Discrete mathematics}, 308(14):3048--3055, 2008.

\bibitem{szollHosi2018constructions}
Ferenc Sz{\"o}ll{\H{o}}si and Patric~RJ {\"O}sterg{\aa}rd.
\newblock Constructions of maximum few-distance sets in euclidean spaces.
\newblock {\em arXiv preprint arXiv:1804.06040}, 2018.

\bibitem{xianglin2012proof}
Wei Xianglin.
\newblock A proof of erd{\H{o}}s-fishburn's conjecture for $ g (6)= 13$.
\newblock {\em The Electronic Journal of Combinatorics}, pages P38--P38, 2012.

\end{thebibliography}
